\documentclass{amsart}

\usepackage{color}

\usepackage{epsfig}
\usepackage{amsmath}
\usepackage{amssymb}
\usepackage{amscd}
\usepackage{graphicx}
\usepackage{color}
\usepackage{tikz}
\usetikzlibrary{arrows,shapes,positioning,shadows,trees}
\tikzset{
  basic/.style  = {draw, text width=2cm, drop shadow, font=\sffamily, rectangle},
  root/.style   = {basic, rounded corners=2pt, thin, align=center,
                   fill=white!30, text width=12em},
  level 2/.style = {basic, rounded corners=6pt, thin,align=center, fill=white!60,
                   text width=9em},
  level 3/.style = {basic, rounded corners=6pt, thin,align=center, fill=white,
                   text width=9em}
}

\newcommand\mL{L\kern-0.08cm\char39}

\newtheorem{thm}{Theorem}[section]
\newtheorem*{thm*}{Weiss-Akin-Glasner Theorem}
\newtheorem{lem}[thm]{Lemma}
\newtheorem{prop}[thm]{Proposition}

\newtheorem{cor}[thm]{Corollary}

\newtheorem{rem}[thm]{Remark}

\newtheorem*{ques A}{Question A}
\newtheorem*{ques B}{Question B}
\newtheorem*{thm A}{Theorem FIP}

\def \N {\mathbb N}

\newcommand{\diam}{{\mathrm{diam}}}

\newcommand{\Tran}{{\mathrm{Tran}}}
\newcommand{\Asym}{{\mathrm{Asym}}}
\newcommand{\orb}{{\mathrm{orb}}}

\newcommand{\Prox}{{\mathrm{Prox}}}

\newcommand{\norm}[1]{\lVert#1\rVert}

\numberwithin{equation}{section}

\makeatletter
\@namedef{subjclassname@2010}{%
  \textup{2010} Mathematics Subject Classification}
\makeatother

\begin{document}

\title{Finite intersection property  {\footnotesize and}   dynamical compactness}

\subjclass[2010]{}

\author{Wen Huang, Danylo Khilko, Sergi{\u\i} Kolyada, \\ Alfred Peris and Guohua Zhang}

\address{Department of Mathematics, Sichuan University,
Chengdu, Sichuan 610064, China}

\address{School of Mathematical Sciences, University of Science and Technology of
China, Hefei, Anhui 230026, China}

\email{wenh@mail.ustc.edu.cn}

\address{Faculty of Mechanics and Mathematics, National Taras Shevchenko University of Kyiv,
Academician Glushkov prospectus 4-b, 03127 Kyiv, Ukraine}

\email{dkhilko@ukr.net}

\address{Institute of Mathematics, NASU, Tereshchenkivs'ka 3, 01601 Kyiv, Ukraine}

\email{skolyada@imath.kiev.ua}

\address{IUMPA, Universitat Polit\`ecnica de Val\`encia, Departament de Matem\`atica Aplicada, Edifici 7A, 46022 Val\`encia, Spain}

\email{aperis@mat.upv.es}

\address{School of Mathematical Sciences and LMNS, Fudan University and Shanghai Center for Mathematical Sciences, Shanghai 200433, China}

\email{chiaths.zhang@gmail.com}

\begin{abstract}
Dynamical compactness with respect to a family as a new concept of chaoticity of a dynamical system was introduced and discussed in \cite{HKhKZh2015}.
In this paper we continue to investigate this notion. In particular, we prove that all dynamical systems are dynamically compact with
respect to a Furstenberg family if and only if this family has the finite intersection property. We investigate weak mixing and
weak disjointness by using the concept of dynamical compactness. We also explore further difference between transitive compactness and weak mixing.
As a byproduct, we show that the $\omega_{\mathcal{F}}$-limit  and the $\omega$-limit sets of a point may have quite different
topological structure. Moreover, the equivalence between multi-sensitivity,
sensitive compactness and transitive sensitivity is established for a minimal system. Finally, these notions are also explored in the context of linear dynamics.
\end{abstract}



\dedicatory{This paper is dedicated to Professor Ethan Akin on the occasion of his 70th birthday.}


\maketitle

\markboth{}{}


\section{Introduction}



By a (\emph{topological}) \emph{dynamical system} $(X,T)$ we mean a compact metric space $X$ with a metric $d$ and a continuous
self-surjection $T$ of $X$. We say it \emph{trivial} if the space is a singleton.
\emph{Throughout this paper, we are only interested in a nontrivial dynamical system, where the state space is a compact metric space
without isolated
points.}


This paper is a continuation of the research carried out in \cite{HKhKZh2015},  where
the authors discuss a dynamical property (called dynamical compactness) and examine it firstly for transitive compactness.
Some results of this paper can be considered as a contribution to \emph{dynamical topology} -- an area of the theory of dynamical systems
in which
the topological properties of maps that can be described in dynamical terms.

Let $\mathbb{Z}_{+}$ be the set of all nonnegative integers and $\mathbb{N}$ the set of all positive integers.
Before going on, let us recall the notion of a Furstenberg family from \cite{Akin1997}.
Denote by $\mathcal{P} =
\mathcal{P}(\mathbb{Z}_{+})$ the set of all subsets of $\mathbb{Z}_{+}$.
A subset $\mathcal{F}\subset \mathcal{P}$ is a (\emph{Furstenberg}) \emph{family}, if it
is \emph{hereditary upward}, that is, $F_1 \subset F_2$ and $F_1 \in \mathcal{F}$ imply $F_2 \in \mathcal{F}$.
Any subset $\mathcal{A}$ of $\mathcal{P}$ clearly generates a family $\{F\in \mathcal{P}: F\supset A\ \text{for some}\ A\in \mathcal{A}\}$.
Denote by $\mathcal{B}$ the family of all infinite
subsets of $\mathbb{Z}_{+}$, and by $\mathcal{P}_+$ the family of all nonempty subsets of $\mathbb{Z}_+$.
For a family $\mathcal{F}$,
the \emph{dual family} of $\mathcal{F}$,
denoted by $k\mathcal{F}$, is defined as $$\{F\in \mathcal{P}: F\cap F'\ne \varnothing \mbox{ for any }
F' \in \mathcal{F}\}.$$
A family $\mathcal{F}$ is \emph{proper} if it is a proper subset of $\mathcal{P}$, that is, $\mathbb{Z}_+\in \mathcal{F}$ and
$\varnothing\notin \mathcal{F}$.
 By a \emph{filter} $\mathcal{F}$ we mean a proper family closed under intersection, that is, $F_1, F_2\in \mathcal{F}$
 implies $F_1\cap F_2\in \mathcal{F}$. A filter is \emph{free} if the intersection of all its elements
 is empty. We extend this concept, a family $\mathcal{F}$ is called \emph{free} if the intersection of all
 elements of $\mathcal{F}$ is empty.

For any $F\in \mathcal{P}$, every point $x\in X$
and each subset $G\subset X$,
 we define $T^{F} x= \{T^i x: i\in F\}$, $n_T(x, G)=\{n\in \mathbb{Z}_{+}: T^n x \in G\}$.
The \emph{$\omega$-limit set of $x$ with respect to $\mathcal{F}$} (see \cite{Akin1997}), or shortly the \emph{ $\omega_\mathcal{F}$-limit
set of $x$}, denoted by $\omega_{\mathcal{F}}(x)$\footnote{Remark that
the notation $\omega_{\mathcal{F}}(x)$ used here is different from the one used in \cite{Akin1997} (the notation $\omega_{\mathcal{F}}(x)$ used here is in fact $\omega_{k \mathcal{F}}(x)$ introduced in \cite{Akin1997}). As this paper is a continuation of the research in \cite{HKhKZh2015}, in order to avoid any confusion of notation or concept,
we will follow the ones used in \cite{HKhKZh2015}.}, is defined as
$$\bigcap_{F \in \mathcal{F}} \overline{T^{F}x}
= \{ z\in X:    n_T(x,G) \in k\mathcal{F} \ \mbox{for every neighborhood} \ G \ \mbox{of} \ z \}.$$
Let us remark that not always $\omega_{\mathcal{F}}(x)$ is a subset of the \emph{$\omega$-limit set} $\omega_{T}(x)$, which is defined as
$$\bigcap_{n=1}^\infty \overline{\{T^kx: k\ge n\}}
= \{ z\in X: N_T(x,G) \in \mathcal{B}  \ \mbox{for every neighborhood} \ G \ \mbox{of} \ z \}. $$ For instance, if each element
of $\mathcal{F}$
contains $0$ then any point $x\in \omega_{\mathcal{F}}(x)$. But, as well known,  a point $x\in \omega_{T}(x)$ if and only if $x$ is a
recurrent point\footnote{A point $x\in X$ is called \emph{recurrent} if $x\in \omega_T(x).$}  of $(X,T)$. Nevertheless, if a family $\mathcal{F}$ is free, then
$\omega_{\mathcal{F}}(x) \subset \omega_{T}(x)$ for any point $x\in X$ and if $(X,T)$ has a nonrecurrent point, then the converse is true (see Proposition \ref{201604161045a}).

 A dynamical system $(X,T)$ is called \emph{compact with respect to $\mathcal{F}$}, or shortly
\emph{dynamically compact}, if the $\omega_\mathcal{F}$-limit set
$\omega_{\mathcal{F}}(x)$ is nonempty for all $x\in X$.

H. Furstenberg started a systematic study of transitive
systems in his paper on disjointness in topological dynamics and ergodic
theory \cite{Furstenberg1967}, and the theory was further developed in \cite{FW1978} and \cite{Furstenberg1981}.
Recall that the system $(X, T)$ is
 (\textit{topologically}) \emph{transitive} if $N_T (U_1, U_2)= \{n\in \mathbb{Z}_+: U_1\cap T^{-n}U_2\neq \varnothing\}\ (= \{n\in \mathbb{Z}_+: T^n U_1\cap U_2\neq \varnothing\})\in \mathcal{P}_+$
for any \emph{opene}\footnote{Because we so often have to refer to
open, nonempty subsets, we will call such subsets \emph{opene}.} subsets $U_1, U_2\subset X$, equivalently, $N_T (U_1, U_2)\in \mathcal{B}$
for any opene subsets $U_1, U_2\subset X$.

In \cite{HKhKZh2015} the authors consider one of possible dynamical compactness --- transitive compactness, and its
relations with well-known chaotic properties of dynamical systems.
 Let $\mathcal{N}_T$ be the set of all subsets of $\mathbb{Z}_+$ containing some $N_T (U, V)$, where $U, V$ are opene subsets of $X$.
A dynamical system $(X, T)$ is called
\emph{transitive compact}, if for any point $x\in X$ the $\omega_{\mathcal{N}_T}$-limit set
$\omega_{\mathcal{N}_T}(x)$ is nonempty,
in other words,
for any point $x\in X$ there exists
a point $z \in X$ such that $$ n_T(x,G) \cap N_T(U,V)\ne \varnothing$$
 for any  neighborhood
$G$ of $z$ and any opene subsets $U,V$ of $X$.

 Let $(X, T)$ and $(Y, S)$ be two dynamical systems and $k\in \mathbb{N}$. The product system $(X\times Y, T\times S)$ is defined naturally,
and denote by $(X^k, T^{(k)})$ the product system of $k$ copies of the system $(X, T)$.
        Recall that the system $(X, T)$ is \emph{minimal} if it does not admit a nonempty, closed, proper subset $K$ of $X$
        with $T K\subset K$,
and is \emph{weakly mixing} if the product system $(X^2, T^{(2)})$ is transitive.
Any transitive compact system is obviously topologically transitive, and observe that each weakly mixing system
is transitive compact (\cite{AK}). In fact, as it was shown in \cite{HKhKZh2015}, each of notions are different in general
and equivalent for minimal systems.

Recall a very useful notion of weakly mixing subsets of a system, which was introduced in \cite{BlH2008} and further discussed in \cite{OZ1} and \cite{OZ2}. The notion of weakly mixing subsets can be regraded as a local version
of weak mixing. Among many very interesting properties let us mention just one of them -- positive topological entropy
of a dynamical system implies the existence of weakly mixing sets (see \cite{LiYe} for details). A nontrivial closed subset $A\subset X$ is called
\emph{weakly mixing} if for every $k\geq 2$
and any opene sets $U_1, \dots ,U_k, V_1, \dots , V_k$ of $X$ with $U_i\cap A \neq \varnothing$, $V_i\cap A \not = \varnothing$,
for any $i=1, \dots ,k$, one has that $\bigcap_{i= 1}^kN_T(U_i\cap A, V_i)\neq \varnothing.$
Let $A$ be a weakly mixing subset of $X$ and let $\mathcal{N}_T(A)$ be the set of all subsets of $\mathbb{Z}_+$ containing some $N_T (U\cap A, V)$, where $U, V$ are opene subsets of $X$ intersecting $A$.

The notion of sensitivity was first used by Ruelle \cite{Ru}, which captures the idea that in a chaotic system a small change in
the initial condition can cause a big change in the trajectory.
According to the works by Guckenheimer \cite{Gu}, Auslander and Yorke \cite{AuslanderYorke} a dynamical
system $(X,T)$ is called
{\it sensitive} if there exists $\delta> 0$ such that for every $x\in X$ and every
neighborhood $U_x$ of $x$, there exist $y\in U_x$ and $n\in \mathbb{N}$ with $d(T^n x,T^n y)> \delta$.
Such a $\delta$ is called a \emph{sensitive constant} of $(X, T)$.
Recently in \cite{Subrahmonian2007} Moothathu initiated a way to measure the sensitivity of a dynamical system, by checking how large is the
set of
nonnegative integers for which the sensitivity occurs (see also \cite{liuheng}).  For a positive $\delta$ and a subset $U\subset X$ define
$$ S_T (U, \delta)= \{n\in \mathbb{Z}_+: \mbox { there are } x_1, x_2\in U\ \mbox{such that}
\ d (T^n x_1, T^n x_2)>
\delta\}.$$
A dynamical system $(X, T)$ is called \emph{multi-sensitive} if there exists  $\delta> 0$ such that
 $\bigcap_{i= 1}^k S_T (U_i, \delta)\neq \varnothing$ for any finite collection of
opene $U_1, \dots, U_k\subset X$. Such a $\delta$ is called a \emph{constant of multi-sensitivity} of $(X, T)$.

\medskip

Recall that a collection $A$ of subsets of a set $Y$ has the \emph{finite intersection property} (FIP) if the intersection
of any finite subcollection of $A$ is nonempty. The FIP is useful in formulating an alternative definition
of compactness of a topological space: a topological space is compact if and only if every collection of closed subsets satisfying the FIP has a nonempty
intersection itself (see, for instance \cite{Edw1995, Kelley}).

Obviously that a filter (say $\mathcal{N}_T$, when $(X,T)$ is weakly mixing), the family $\mathcal {N}_T(A)$ for
a weakly mixing subset $A$ of  $(X,T)$  and the family $\mathcal {S}_T(\delta)$ when $(X, T)$ is a multi-sensitive system (with a
constant of multi-sensitivity $\delta> 0$) have FIP.  Since all of these families are also free, actually they  have the
\emph{strong finite intersection property}
(SFIP), i.e., if the intersection over any finite subcollection
of the family is infinite (see Proposition 2.2).

In fact we can say
more --- the FIP is useful in characterizing the dynamical compactness (see Theorem \ref{main}).

\begin{thm A}
All dynamical systems are dynamically compact with respect to $\mathcal{F}$ if and only if the family
$\mathcal{F}$ has the finite intersection property.
\end{thm A}

We also introduce  two new stronger versions of sensitivity: sensitive compactness and transitive sensitivity.
 Denote by $\mathcal{S}_T(\delta)$ the set of all subsets of $\mathbb{Z}_+$ containing $S_T (U, \delta)$ for some  $\delta> 0$ and
 opene $U\subset X$.
 We will call the
system $(X, T)$ \emph{transitively sensitive} if there exists  $\delta> 0$ such that  $S_T (W, \delta)\cap N_T (U, V)\neq \varnothing$
for any opene subsets $U, V, W$ of $X$; and
\emph{sensitive compact}, if there exists $\delta> 0$ such that for any point $x\in X$ the $\omega_{\mathcal{S}_T(\delta)}$-limit set
$\omega_{\mathcal{S}_T(\delta)}(x)$ is nonempty,
in other words,
for any point $x\in X$ there exists
a point $z\in X$ such that
 $$ n_T(x,G)\cap S_T (U, \delta)\neq \varnothing$$
 for any  neighborhood
$G$ of $z$ and any opene $U$ of $X$.

\medskip

The paper is organized as follows. In Section 2 we recall some basic concepts and properties used in later
discussions from topological dynamics. In Section 3
we obtain some general results concerning dynamical compactness. In particular we show that all dynamical systems are dynamically compact
with respect to a Furstenberg family if and only if this family has the finite intersection property
(Theorem \ref{main}).

In Section 4 we discuss two stronger versions of sensitivity: transitive sensitivity and sensitive compactness.
It was shown that each weakly mixing system is transitively sensitive (Proposition \ref{201604181936}), and in fact we can characterize
transitive sensitivity of a general dynamical system in terms of dynamical compactness (Proposition \ref{201604092152}). Furthermore,
all of the multi-sensitivity, sensitive compactness and transitive sensitivity are equivalent for a minimal system
(Theorem \ref{201604092203}).
Even though each minimal transitive compact system is multi-sensitive,
nevertheless, there are many minimal multi-sensitive systems which are not transitive compact.
Observe that the sensitivity of a dynamical system can be lifted up from a factor to an extension by an almost open
factor map between transitive systems by \cite[Corollary 1.7]{GlasnerWeiss1993}. We prove that the transitive sensitivity
can be lifted up to an extension from a factor by an almost one-to-one factor map and that the transitive sensitivity is
projected from an extension to the sensitivity of a factor by a weakly almost one-to-one factor map (Lemma \ref{201604101854}).

In Section 5 we show that dynamical compactness can be used to characterize the weak disjointness of dynamical
systems (Theorem \ref{th0201}). We also improve the result of Jian Li \cite{Li}: weak mixing implies $\mathcal{F}_{ip}$-point
transitivity in terms of transitive compactness (Proposition \ref{0203}).

In Section 6 the further difference between weak mixing and transitive compactness is explored.
 Precisely, there is a totally transitive, non weakly mixing, transitive compact
system (Theorem \ref{th0601}); and in fact any compact metric space can be realized as the $\omega_{\mathcal{N}_T}$-limit set of a
non totally
transitive, transitive compact system $(X, T)$ (Theorem \ref{20151026}). As a byproduct, we show that the $\omega_{\mathcal{N}_T}$-limit
sets and the $\omega$-limit sets have quite different topological structures for a general dynamical compact system $(X,T)$. At the end
of this section we add one more chaotical property of transitive compact systems (in additional to already known
from \cite{HKhKZh2015}): transitive compactness implies Li-Yorke chaos (Proposition \ref{0200a}).

In Section 7 we consider the dynamics of linear operators on infinite dimensional spaces in relation to the properties studied
in previous sections. In particular, we show the equivalence of the topological weak mixing property with a weak version of transitive
compactness (Theorem
\ref{wmoptc}). Some results on sensitivity are also obtained.

\vskip 10pt

\noindent {\bf Acknowledgements.}  The first and third authors acknowledge the hospitality of the School of Mathematical Sciences of
the Fudan University, Shanghai.
The third author also acknowledges the hospitality of the Departament de Matem\`atica Aplicada of the
Universitat Polit\`ecnica de Val\`encia and the Department of Mathematics of the Chinese University of Hong Kong.

The first author was supported by NNSF of China (11225105, 11431012); the fourth author was supported by MINECO, Project MTM2013-47093-P,
and by GVA, Project
PROMETEOII/2013/013; and the fifth author
was supported by NNSF of China (11271078).

\section{Preliminaries}

In this section we recall standard concepts and results used in later discussions.

\subsection{Basic concepts in topological dynamics}

 Recall that $x\in X$ is a \emph{fixed point} if $T x= x$, and an \emph{$\mathcal{F}$-transitive point of $(X, T)$} \cite{Li}
 if $n_T(x, U)\in \mathcal{F}$ for any opene subset $U$ of $X$.
 It is a trivial observation that if a family $\mathcal{F}$ admits an
$\mathcal{F}$-transitive
dynamical system $(X,T)$ without isolated points,
then $\mathcal{F}$ is free.
Since $k(k\mathcal{F})= \mathcal{F}$, it is easy to
see that $x\in X$ is an $\mathcal{F}$-transitive point of $(X, T)$ if and only if $\omega_{k\mathcal{F}}(x)=X$.
 Denote by $\Tran_\mathcal{F} (X, T)$ the set of all $\mathcal{F}$-transitive points of $(X, T)$.
The system
$(X, T)$ is \emph{$\mathcal{F}$-point transitive} if $\Tran_\mathcal{F} (X, T)\neq \varnothing$, and is \emph{$\mathcal{F}$-transitive} if $N_T(U,V)\in \mathcal{F}$ for any opene subsets $U,V$ of $X$.
Write $\Tran (X, T)= \Tran_{\mathcal{P}_+} (X, T)$ for short, and we also call the point $x$ \emph{transitive} if $x\in \Tran (X, T)$, equivalently, its \emph{orbit} $\orb_T (x)= \{T^n x: x= 0, 1, 2, \dots\}$ is dense in $X$.
Since $T$ is surjective, the system $(X, T)$ is transitive if and only if $\Tran (X, T)$ is a dense
    $G_\delta$ subset of $X$.

    In general, a subset $A$ of $X$ is \emph{$T$-invariant} if $TA = A$, and \emph{positively $T$-invariant} if $T A\subset A$.
    If $A$ is a closed, nonempty, $T$-invariant subset then $(A,T|_A)$ is called the associated \emph{subsystem}.
    A \emph{minimal subset} of $X$ is a closed, nonempty, $T$-invariant subset such that the associated subsystem is minimal.
    Clearly, $(X,T)$ is minimal
if and only if $\Tran(X,T) = X$,
     if and only if it admits no a proper, closed, nonempty, positively $T$-invariant subset.
     A point $x \in X$ is called \emph{minimal} if it lies in some minimal subset. In this case, in order to emphasize the underlying system $(X, T)$ we also say that $x\in X$ is a \emph{minimal point of $(X, T)$}.
     Zorn's Lemma implies that every closed, nonempty, positively $T$-invariant set contains a minimal set.

A pair of points $x, y \in X$ is called \emph{proximal} if $\liminf_{n\rightarrow \infty} d (T^n x, T^n y)= 0$.
In this case each of points from the pair is said to be \emph{proximal} to another. Denote by $\Prox _T (X)$ the set of all proximal pairs of points.
For each $x\in X$, denote by $\Prox _T (x)$, called the \emph{proximal cell of $x$}, the set of all points
which are proximal to $x$.
Recall that a dynamical system $(X,T)$ is called \emph{proximal}
if $\Prox_T (X)= X\times X$. The system $(X,T)$ is
proximal if and only if $(X,T)$ has the unique fixed point, which is the only minimal point of $(X,T)$ (e.g. see \cite{AK}).

The opposition to the notion of sensitivity is the concept of equicontinuity.
Recall that $x \in X$ is an \emph{equicontinuity point} of $(X, T)$ if for every $\varepsilon > 0$ there exists a $\delta > 0$ such that $d (x, x') < \delta$
implies $d (T^n x, T^n x')< \varepsilon$ for any $n\in \mathbb{Z}_+$.
Denote by $\text{Eq} (X, T)$ the set of all equicontinuity points of $(X, T)$.
The system $(X, T)$ is called \emph{equicontinuous} if $\text{Eq} (X, T)= X$.
 Each dynamical system admits a maximal equicontinuous factor.
       Recall that by a \emph{factor map} $\pi: (X, T)\rightarrow (Y, S)$
       between dynamical systems  $(X, T)$ and $(Y, S)$,
        we mean that $\pi: X\rightarrow Y$ is a continuous
 surjection with $\pi\circ T= S\circ \pi$. In this case, we call $\pi: (X, T)\rightarrow (Y, S)$
    an \emph{extension}; and $(X, T)$  an \emph{extension} of
  $(Y, S)$, $(Y, S)$ a \emph{factor} of $(X, T)$.

\subsection{Basic concepts of Furstenberg families}

In this subsection we recall from \cite{Akin1997} basic concepts about Furstenberg families.

Let $F\in \mathcal{P}$. Recall that a subset $F$ is \emph{thick} if it contains arbitrarily long runs of positive integers. Denote by $\mathcal{F}_{\text{thick}}$ the set of all thick subsets of $\mathbb{Z}_+$, and define $\mathcal{F}_{\text{syn}}= k \mathcal{F}_{\text{thick}}$.
Each element of $\mathcal{F}_{\text{syn}}$ is said to be \emph{syndetic}, equivalently, $F$ is syndetic if and only if there is $N\in \mathbb{N}$ such that $\{i, i + 1, \dots , i + N \} \cap F \not = \varnothing$
for every $i \in \mathbb{Z}_{+}$. We say that $F$ is \emph{thickly syndetic} if for every $N\in \mathbb{N}$ the positions where
length $N$ runs begin form a syndetic set. Denote by $\mathcal{F}_{\text{cof}}$ the set of all cofinite subsets of $\mathbb{Z}_+$.
Note that by the classic result of Gottschalk a point
 $x\in X$ is minimal if and only if $n_T(x, U)= \{n\in \mathbb{Z}_+: T^n x\in U\}$ is syndetic for any neighborhood $U$ of $x$. Hence, for any minimal system $(X, T)$, the subset $N_T (U, V)$ is syndetic for any opene subsets $U, V$ of $X$.

    Recall that a family $\mathcal{F}$ is proper if it is a proper subset of $\mathcal{P}$, that is, $\mathbb{Z}_+\in \mathcal{F}$ and $\varnothing\notin \mathcal{F}$.
 By a filter $\mathcal{F}$ we mean a proper family closed under intersection, that is, $F_1, F_2\in \mathcal{F}$ implies $F_1\cap F_2\in \mathcal{F}$.
For families $\mathcal{F}_1$ and
$\mathcal{F}_2$, we define the family $\mathcal{F}_1\cdot \mathcal{F}_2:= \{F_1\cap F_2: F_1\in \mathcal{F}_1, F_2\in \mathcal{F}_2\}$ and call it the \emph{interaction of $\mathcal{F}_1$ and $\mathcal{F}_2$}. Thus we have
$\mathcal{F}_1 \cup \mathcal{F}_2 \subset \mathcal{F}_1 \cdot \mathcal{F}_2$; and it is easy to check that $\mathcal{F}$
is a filter if and only if $\mathcal{F}= \mathcal{F}\cdot \mathcal{F}$, and $\mathcal{F}_1\cdot \mathcal{F}_2$  is
proper if and only if $\mathcal{F}_2 \subset k\mathcal{F}_1$.

 For each $i\in \mathbb{Z}_+$, define $g^i: \mathbb{Z}_+\rightarrow \mathbb{Z}_+, j\mapsto i+ j$. Recall that a family $\mathcal{F}$ is
 \emph{$+$ invariant} if for every $i\in \mathbb{Z}_+$, $F\in \mathcal{F}$ implies $g^{i} (F)\in \mathcal{F}$; \emph{$-$ invariant} if
 for every $i\in \mathbb{Z}_+$, $F\in \mathcal{F}$ implies $g^{- i} (F)\in \mathcal{F}$, where $g^{-i}(F)=(g^i)^{-1}(F)=
 \{ j-i:j\in F,j\ge i\}$; and
  \emph{translation invariant} if it is both $+$ and $-$ invariant, equivalently, for every $i\in \mathbb{Z}_+$, $F\in \mathcal{F}$ if and only if $g^{- i} (F)\in \mathcal{F}$.

  As $g^{- i} (g^i A)= A$ and $g^i (g^{- i} A)\subset A$ for any $i\in \mathbb{Z}_{+}$, it is easy to obtain that the family $\mathcal{F}$ is $+$ ($-$, translation, respectively) invariant if and only if $k \mathcal{F}$ is $-$ ($+$, translation, respectively) invariant (see for example \cite[Proposition 2.5.b]{Akin1997}).
  And then we have:

\begin{prop} \label{201604161045}
Let $x\in X$. Then
$T \omega_\mathcal{F} (x)\subset \omega_\mathcal{F} (T x)$. Additionally, if $\mathcal{F}$ is $-$ ($+$, translation, respectively) invariant then $\omega_\mathcal{F} (T x)\subset$ ($\supset$, $=$, respectively) $\omega_\mathcal{F} (x)$.
\end{prop}

\begin{proof}
Since the other items are
alternative versions of \cite[Proposition 3.6]{Akin1997} in our notations, it suffices to prove that if $\mathcal{F}$ is $+$ invariant then $\omega_\mathcal{F} (T x)\supset \omega_\mathcal{F} (x)$.

For each $y\in \omega_\mathcal{F} (x)$ take an arbitrary neighborhood $U$ of $y$, and let $F\in \mathcal{F}$. Then $g^{1} (F)= \{i+ 1\in \mathbb{Z}_+: i\in F\}\in \mathcal{F}$ as $\mathcal{F}$ is $+$ invariant, and hence $n_T (x, U)\cap g^{1} (F)\neq \varnothing$, thus $\varnothing\neq g^{- 1} (n_T (x, U)\cap g^{1} (F))= n_T (T x, U)\cap F$. It follows $y\in \omega_\mathcal{F} (T x)$ from the arbitrariness of $U$ and $F$, which finishes the proof.
\end{proof}

\begin{prop} \label{201604161045a}
Let $(X,T)$ be a dynamical system and let $\mathcal{F}$ be a family.
\begin{enumerate} \item[(i)] If $\mathcal{F}$ is free, then
$\omega_{\mathcal{F}}(x) \subset \omega_{T}(x)$ for any $x\in X$. Moreover, if $(X,T)$ has a nonrecurrent point, then the converse is
true.
\item[(ii)]
If $\mathcal{F}$ is free and has FIP then it has SFIP.
\end{enumerate}

\end{prop}

\begin{proof}

(i) Suppose $\omega_{\mathcal{F}}(x)\neq \varnothing$ and take a point $y \in
\omega_{\mathcal{F}} (x):=\bigcap_{F \in \mathcal{F}} \overline{T^{F}x}$. Let us show that $y \in
\omega_{T} (x)$.  Since $\mathcal{F}$ is free, $\bigcap_{F \in \mathcal{F}} T^{F}x=\varnothing$. Otherwise
there is $m \in \mathbb{Z}_+$ such that $T^mx \in T^Fx$ for all $F \in \mathcal{F}$, in other words
$m \in \bigcap_{F \in \mathcal{F}}F$, a contradiction. It means that  there is  $F \in \mathcal{F}$
such that $y \notin T^Fx$.  So $y \in \overline{T^{F}x} \setminus T^Fx$, and $F$ is
infinite. Therefore there exists an infinite sequence $T^{n_1}x, \dots, T^{n_i}x,\dots, $ which converges
to $y$. Hence $y \in \omega_T(x)$.

Now, let $(X,T)$ be a dynamical system with a nonrecurrent point $x_0\in X$, let $\mathcal{F}$ be a family and
$\omega_{\mathcal{F}}(x) \subset \omega_{T}(x)$ for any $x\in X$. Suppose $\mathcal{F}$ is not free.
It means there is a $k\in \mathbb{Z}_{+}$ that lies in each element of $\mathcal{F}$. Then obviously that $T^k(x)\in
\omega_{\mathcal{F}}(x)$ for any $x\in X$. Since $x_0$ is nonrecurrent, $x_0\notin \omega_T(x_0)\ne X$. It is well known that $\omega_T(x)$ is
$T$-invariant, therefore $\omega_T(y)=\omega_T(x_0)$ for any  $y\in \{T^{-i}(x_0): i \in \mathbb{Z}_{+}\}$, and
$y\notin \omega_T(y)$.  Take a point $x_k\in X$ with $T^k(x_k)=x_0$. As we know $x_0 \in \omega_{\mathcal{F}}(x_k)$. But $\omega_{\mathcal{F}}(x_k)\subset \omega_T(x_k)=  \omega_T(x_0)$, a contradiction.

(ii) Suppose that for some $F_1, \dots F_k \in \mathcal{F}$ $\bigcap_{i=1}^k F_i=\{n_1, n_2, \dots n_m\}$.
Since $\mathcal{F}$ is free for each $k \in \mathbb{Z}_+$ there is $G_k \in \mathcal{F}$ such that $k \notin G_k$.
Then $\bigcap_{i=1}^k F_i \cap \bigcap_{j=1}^m G_{n_j}=\varnothing$, contradiction.
\end{proof}

\subsection{The concept of an almost one-to-one map}

Let $\phi: X\rightarrow Y$ be a continuous surjective map from a compact metric space $X$ onto a compact Hausdorff space $Y$.
Recall that $\phi$ is
\emph{almost open} if $\phi (U)$ has a nonempty interior in $Y$ for any opene $U\subset X$.
Note that each factor map between
minimal systems is almost open \cite[Theorem 1.15]{Auslander1988}, in particular, for a minimal system $(X,T)$ the map $T:X\to X$ is almost open \cite{KST2001}.
Denote by $Y_0\subset Y$ the set of all points $y\in Y$ whose fiber is a singleton.
 Then $Y_0$ is a $G_\delta$ subset of $Y$,
 because
$$Y_0= \{y\in Y: \phi^{- 1} (y)\ \text{is a singleton}\}= \bigcap_{n\in \mathbb{N}} \left\{y\in Y: \diam (\phi^{- 1} (y))<
\frac{1}{n}\right\}$$
and the map $y\mapsto \diam (\phi^{- 1} (y))$ is upper semi-continuous. Here, we denote by $\diam (A)$ the diameter of a subset $A\subset X$.
 Recall that
the function $f: Y\rightarrow \mathbb{R}_+$ is
\emph{upper semi-continuous} if $\limsup\limits_{y\rightarrow y_0} f (y)\le f (y_0)$ for each $y_0\in Y$.
Denote by $X_0\subset X$ the set of all points $x\in X$ such that the pre-image of $\phi (x)$ is a singleton. Then $X_0= \pi^{- 1} (Y_0)$ is a $G_\delta$ subset of $X$.

We call $\phi$ \emph{weakly almost one-to-one} if $Y_0$ is dense in $Y$, and \emph{almost one-to-one}\footnote{Here we use the concept of \emph{almost one-to-one} following \cite{AkinGlasner2001}, and the concept of \emph{almost one-to-one} used in \cite{Downarowicz2005, HKZ, KST2001} is in fact our weakly almost one-to-one.} if $X_0$ is dense in $X$.
It is not hard to show that:
if $\phi$ is weakly almost one-to-one, then for any $\delta> 0$ and any opene subset $U$ of $Y$ there exists opene $V\subset U$ with $\diam (\phi^{- 1} V)< \delta$; and if $\phi$ is almost one-to-one, then for any opene subset $U^*$ of $X$ there exists an opene subset $V^*$ of $Y$ with $\phi^{- 1} V^*\subset U^*$.
Clearly almost one-to-one is much stronger than weakly almost one-to-one. For example, let $X$ be the closed unit interval, define $T (x)= 2 x$ for $x\in [0, \frac{1}{2}]$ and $T (x)= 1$ for $x\in [\frac{1}{2}, 1]$, and then $T: X\rightarrow X$ is clearly not almost one-to-one but weakly almost one-to-one.

For each minimal system $(X, T)$, the map $T: X\rightarrow X$ is weakly almost one-to-one \cite[Theorem 2.7]{KST2001}, and in fact almost one-to-one \cite[Proposition 2.3]{HKZ}.
The following result characterizes the relationship between weakly almost one-to-one and almost one-to-one, which extends \cite[Proposition 2.3]{HKZ}.

\begin{prop} \label{201604081725}
Let $\phi: X\rightarrow Y$ be a continuous surjective map from a compact metric space $X$ onto a compact Hausdorff space $Y$.
 Then $\phi$ is almost one-to-one if and only if it is not only almost open but also weakly almost one-to-one.
\end{prop}

\begin{proof}
Firstly assume that $\phi$ is almost one-to-one. Let $U\subset X$ be an arbitrary opene subset. And then we can take $x_0\in U$ such that the pre-image of $\phi (x_0)$ is a singleton. From this it is easy to see that $\phi (x_0)$ is contained in the interior of $\phi (U)$. This implies that $\phi$ is almost open. The map $\phi$ is clearly weakly almost one-to-one.

Now assume that $\phi$ is not only almost open but also weakly almost one-to-one.
 Let $U\subset X$ be an arbitrary opene subset.
 Since $\phi$ is almost open, $\phi (U)$ has a nonempty interior in $Y$, and then $\phi^{- 1} (y_0)$ is a singleton for some $y_0\in \phi (U)$, as $\phi$ is weakly almost one-to-one. This shows $U\cap X_0\neq \varnothing$, which finishes the proof.
\end{proof}

As a direct corollary, we have:

\begin{cor} \label{201604102314}
Let $\phi: X\rightarrow Y$ and $\pi: Y\rightarrow Z$ be continuous surjective maps between compact metric spaces. Then the composition map $\pi\circ \phi: X\rightarrow Z$ is almost one-to-one if and only if both $\phi$ and $\pi$ are almost one-to-one.
\end{cor}

\begin{proof}
Denote by $X_0$ ($X_1$, respectively) the set of all points $x\in X$ such that the pre-image of $(\pi\circ \phi) (x)$ ($\phi (x)$, respectively) is a singleton. Denote by $Z_0$ ($Z_1$, respectively) the set of all points $z\in Z$ whose $\pi\circ \phi$-fibers ($\pi$-fibers, respectively) are singletons. All of them are $G_\delta$ subsets. Moreover, $X_0= X_1\cap \phi^{- 1} (\pi^{- 1} Z_1)$. In fact, $x\in X_0$ if and only if $\{x\}= (\pi\circ \phi)^{- 1} (\pi\circ \phi (x))= \phi^{- 1} (\pi^{- 1} (\pi (\phi x)))$, if and only if $\pi^{- 1} (\pi (\phi x))= \{\phi (x)\}$ and $\phi^{- 1} (\phi x)= \{x\}$, if and only if $\pi (\phi x)\in Z_1$ and $x\in X_1$.

First assume that $\pi\circ \phi$ is almost one-to-one, and then by Proposition \ref{201604081725}: $X_0$ is a dense subset of $X$, $Z_0$ is a dense subset of $Z$ and the map $\pi\circ \phi$ is almost open. Note that $X_0\subset X_1$ and $Z_0\subset Z_1$, we have that $X_1$ is dense in $X$ and $Z_1$ is dense in $Z$. Hence $\phi$ is almost one-to-one. Furthermore, as the map $\pi\circ \phi$ is almost open, for any opene $V\subset Y$ one has that $\pi (V)= (\pi\circ \phi) (\phi^{- 1} V)$ has a nonempty interior in $Z$, which implies that $\pi$ is almost one-to-one by Proposition \ref{201604081725}.

Now assume that both $\phi$ and $\pi$ are almost one-to-one. Then $X_1$ is a dense $G_\delta$ subset of $X$ and $Z_1$ is a dense $G_\delta$ subset of $Z$. Moreover, by Proposition \ref{201604081725} both $\phi$ and $\pi$ are almost open, and then the continuous surjection $\pi\circ \phi$ is also almost open, which implies that $(\pi\circ \phi)^{- 1} (Z_1)$ is also a dense $G_\delta$ subset of $X$. Thus, $X_0= X_1\cap \phi^{- 1} (\pi^{- 1} Z_1)$ is a dense $G_\delta$ subset of $X$, that is, the composition map $\pi\circ \phi: X\rightarrow Z$ is almost one-to-one. This finishes the proof.
\end{proof}

Let $\pi: (X, T)\rightarrow (Y, S)$ be a factor map between dynamical systems.
If the map $\pi: X\rightarrow Y$ is almost one-to-one (weakly almost one-to-one, respectively),
then we also call $(X, T)$ an \emph{almost one-to-one extension} (a \emph{weakly almost one-to-one extension}, respectively) \emph{of $(Y, S)$}.
The main result of \cite{HKZ} states that a minimal system is either multi-sensitive or a weakly almost one-to-one extension of its maximal equicontinuous factor. This is an analog of the well-known Auslander-Yorke dichotomy theorem: a minimal system is either sensitive or equicontinuous.

\subsection{Symbolic dynamics}

Let $A$ be a nonempty finite set. We call $A$ the \emph{alphabet}
and elements of $A$ are \emph{symbols}. The \emph{full \emph{(}one-sided\emph{)} $A$-shift} is defined as
$$\Sigma= \{x = \{x_i\}_{i=0}^\infty: x_i \in A \text{ for all } i \in \mathbb{Z}_+\},$$
where
we equip $A$ with the discrete topology and $\Sigma$ with the product topology, and the \emph{shift
map} $\sigma : \Sigma \to \Sigma$  is a continuous surjection given by
$$ x = \{x_i\}_{i=0}^\infty \mapsto \sigma x = \{x_{i+1}\}_{i=0}^\infty,$$
that is, $\sigma (x)$ is the sequence obtained by dropping the first symbol of $x$.
 Usually we write an element of $\Sigma$ as $x = \{x_i\}_{i=0}^\infty = x_0 x_1 x_2 x_3 \dots$

A \emph{block} $w$ over $\Sigma$ is a finite sequence of symbols and its \emph{length}
is the number of its symbols (denoted by $|w|$). An \emph{$n$-block} stands for a block of length $n$. In general we are only interested in a block $w$ with $|w|\ge 1$ if without any special statement, and denote by $\Sigma^*$ the
set of all blocks over $\Sigma$. The block $w$ is a \emph{subblock} of a block $v= v_1 \dots v_m$ with $v_1, \dots, v_m\in A$ if there exists $1\le i\le j\le m$ with $w= v_i \dots v_j$. The concatenation of two
blocks $u = a_1 \dots a_k$ and $v = b_1 \dots b_l$ is the block $uv = a_1 \dots a_k b_1 \dots b_l$. We write
$u^n$ for the concatenation of $n\geq 1$ copies of a block $u$ and $u^\infty$ for the sequence $uuu \dots \in \Sigma.$ By $x_{[i,j]}$ we denote the block $x_i x_{i+1} \dots x_j$, where $0 \leq i \leq j$ and $x =
\{x_k\}_{k=0}^\infty \in \Sigma$. The subset $X\subset \Sigma$ is called a \emph{subshift} if it
is a closed, nonempty, $\sigma$-invariant subset of $\Sigma$. A \emph{cylinder} of an $n$-block $w \in \Sigma^*$
in a subshift $X$ is the set $C[w] = \{x \in X: x_{[0, n-1]} = w\}$. The collection of all cylinders forms a basis of the topology of
$X$.

\section{Dynamical compactness with respect to an arbitrary family}

Recall that a family  $\mathcal{F}$ has the finite intersection property (FIP) if the intersection
of any finite subcollection of $\mathcal{F}$ is nonempty. The following theorem shows that the FIP is
useful in characterizing the dynamical compactness.

\begin{thm} \label{main}
All dynamical systems are dynamically compact with respect to the family $\mathcal{F}$ if and only if
$\mathcal{F}$ has the finite intersection property.
\end{thm}

\begin{proof}
\textit{Sufficiency}.  Suppose that $\mathcal{F}$ has FIP. Take arbitrary dynamical system $(X, T)$ and let $x\in X$. Obviously the family
$\{\overline{T^{F}x}: F\in \mathcal{F}\}$ also has FIP, and then by compactness of the space $X$ the family has a nonempty intersection itself, i.e., $\omega_{\mathcal{F}}(x)=
\bigcap_{F \in \mathcal{F}} \overline{T^{F}x} \neq \varnothing$. Thus $(X, T)$ is dynamically compact with respect to $\mathcal{F}$.

\textit{Necessity}.
Suppose that the family $\mathcal{F}$ has no FIP. And then there is a collection
$\{F_1, \dots, F_k \} \subset \mathcal{F}$ with $\bigcap_{i=1}^k F_i=\varnothing.$
Let $A=\{a_1, \dots, a_k\}$ be an alphabet and let $(X, T):=(\Sigma, \sigma)$ be the full (one-sided) $A$-shift.
We are going to define a point $x \in X$ with $\omega_\mathcal{F}(x)= \varnothing$.  Let $x_0=a_1$. For any $n \geq 1$ there is $i$ with $n \notin F_i$, else the intersection of $F_1, \dots, F_k$ would
be nonempty. Then define $x_n:=a_i$. Finally, let $x=x_0 x_1 x_2 x_3\dots$ and the construction is finished.

Assume the contrary that we can take $z \in \omega_\mathcal{F}(x)$, and that $z$ begins with $a_i \in A$. Take $G_z=C[a_i]$. As $z \in \omega_\mathcal{F}(x)$ we have $n_T(x, G_z) \cap F_i \neq \varnothing$. But
if $n \in n_T(x, G_z)$,  then $x_n=a_i$ and so $n \notin F_i$ by the construction, a contradiction.
\end{proof}

As we have mentioned in Introduction, obviously that a filter $\mathcal{F}$ (in particular $\mathcal{N}_T$, when $(X,T)$ is weakly mixing), the family $\mathcal {N}_T(A)$ for
a weakly mixing subset $A$ of  $(X,T)$  and the family $\mathcal {S}_T(\delta)$ when $(X, T)$ is a multi-sensitive system (with a
constant of multi-sensitivity $\delta> 0$) have FIP.

Let $\mathcal{F}$ has the finite intersection property. Then there exists an ultrafilter $\mathcal{U}$ (in $\mathcal{P}$) such that
$\mathcal{F}\subset \mathcal{U}$. This result is known as Ultrafilter Lemma (see details and proof in \cite{Halp}). Recall that an
ultrafilter is maximal among all proper filters. As a consequence of this fact we have a natural open question:

\begin{ques A}
 Let $(X,T)$ be a dynamically compact system with respect to a family $\mathcal{F}$ and $\mathcal{F}$ has FIP.
 When $\mathcal{F}$ is a filter, or at least contains a nontrivial
 filter?
\end{ques A}

Especially we address this question to the family  $\mathcal {S}_T(\delta)$. More precisely, when a system $(X,T)$
is dynamically compact with respect to the family $\mathcal{N}_T$ and  $\mathcal{N}_T$  has FIP, then, as well known,
the systems is weakly mixing and $\mathcal{N}_T$ is a filter. Now, let a system $(X,T)$ is dynamically compact with respect
to the family $\mathcal {S}_T(\delta)$ for some $\delta >0$ and $\mathcal {S}_T(\delta)$ has FIP, then the systems is multi-sensitive.
But, the following question is still open -- when is $\mathcal {S}_T(\delta)$ a filter?

A collection $\mathcal{H}\subset \mathcal{F}$ will be called a \textit{base} for $\mathcal{F}$ if for any $F\in \mathcal{F}$ there is $H \in \mathcal{H}$ with $H \subset F$.
We are interested in those families which have \textit{a countable base}, that is, there exists a
base $\mathcal{H}$ which is countable.

Remark that not any Furstenberg family $\mathcal{F}$ has a countable base, for example,
the family $\mathcal{B}$. Assume the contrary that $\mathcal{B}$ admits a countable base $\{F_n: n\in \mathbb{N}\}$. We take $k_1\in F_1$,
and once $k_m\in F_m, m\in \mathbb{N}$ is defined we choose $k_{m+ 1}\in F_{m+ 1}$ with $k_{m+ 1}> k_m+ m+ 1$.
Set $E= \{k_n: n\in \mathbb{N}\}$ and $F= \mathbb{Z}_+\setminus E$.
Then $E\cap F_n\neq \varnothing$ for all $n\in \mathbb{N}$, and $F\supset \{k_m+ m: m\in \mathbb{N}\}$ and hence $F\in \mathcal{B}$,
in particular, there exists no $n\in \mathbb{N}$ with $F_n\subset F$, a contradiction.

Not so hard to show even the existence of a family with FIP, but without a countable base. Nevertheless the families $\mathcal{N}_T$ and $\mathcal{S}_T (\delta)$ have countable bases.
Indeed, we can consider a countable base $\mathcal{U}$ of open sets for the space $X$.
Note that $U_1 \subset U$, $V_1 \subset V$ implies $N_T(U_1, V_1) \subset N_T(U, V)$ and $S_T(U_1,\delta) \subset S_T(U,\delta)$.
Then $\{N_T(U, V): ~U, V \in \mathcal{U}\}$ and $\{S_T(U,\delta):~U \in \mathcal{U}\}$ are countable bases
for $\mathcal{N}_T$ and $\mathcal{S}_T (\delta)$, respectively.

The following is a general result that will be especially useful for families with countable bases.

\begin{prop} \label{orbtranspoint}
Let $(X,T)$ be a dynamical system and let $\mathcal{F}$ be a family such that there exists $x\in \Tran_{\mathcal{F}} (X, T)$. Then $\orb_T(x)\subset \Tran_{\mathcal{F}} (X, T)$.
\end{prop}

\begin{proof} By assumption, given an arbitrary opene $U\subset X$ we have that $n_T(x,U)\in \mathcal{F}$. Thus, for any $m\in\N$,
$$
n_T(T^mx,U)=n_T(x,T^{-m}(U))\in \mathcal{F},
$$
and we conclude that $T^mx\in \Tran_{\mathcal{F}} (X, T)$.
\end{proof}

\begin{prop}  \label{prop0201}
Assume that $\mathcal{F}$ admits a countable base $\mathcal{H}$. Then $\Tran_{k\mathcal{F}} (X, T)$ is a $G_\delta$ subset of $X$. Moreover, the  following are equivalent:
\begin{enumerate}
\item[(1)] The system $(X,T)$ is $k \mathcal{F}$-transitive,
\item[(2)] $\Tran_{k\mathcal{F}} (X, T)$ is a dense $G_\delta$ subset of $X$,
\item[(3)] $\Tran_{k\mathcal{F}} (X, T)\neq \varnothing$.
\end{enumerate}
\end{prop}

\begin{proof}
Let $\mathcal{U}$ be a countable base of the family of all opene subsets of $X$. Then the class $\mathcal{U} \times \mathcal{H}$ is countable, and we enumerate it as  $\{(U_i, F_i): i \in \mathbb{N}\}$. Denote by $T^{- F} U= \bigcup_{n\in F} T^{- n} U$ for any $F\subset \mathbb{Z}_+$ and each $U\subset X$.
Then it is easy to obtain
\begin{equation} \label{201604161338}
\Tran_{k\mathcal{F}}(X,T) = \bigcap_{i=1}^\infty T^{-F_i} U_i.
\end{equation}
In fact, given arbitrary point $x\in X$, $x\in \Tran_{k\mathcal{F}}(X,T)$ if and only if $N_T (x, U)\in k \mathcal{F}$ for any opene subset $U$ of $X$, if and only if $N_T (x, U)\cap F\neq \varnothing$ for any opene subset $U$ of $X$ and each $F\in \mathcal{F}$, if and only if $N_T (x, U_i)\cap F_i\neq \varnothing$ for each $i\in \mathbb{N}$ by the construction. In particular, $\Tran_{k\mathcal{F}} (X, T)$ is a $G_\delta$ subset of $X$.

Thus $(X, T)$ is $k \mathcal{F}$-transitive, if and only if for any $F\in \mathcal{F}$ and arbitrary opene subsets $U, V$ of $X$ we have $N_T (V, U)\cap F\neq \varnothing$ and equivalently $T^{- F} U\cap V\neq \varnothing$, if and only if
    $T^{- F} U$ is an opene dense subset of $X$ for any $F\in \mathcal{F}$ and each opene subset $U$ of $X$, if and only if $\Tran_{k\mathcal{F}} (X, T)$ is a dense $G_\delta$ subset of $X$ by \eqref{201604161338}.

Now we assume $\Tran_{k\mathcal{F}} (X, T)\neq \varnothing$.    Let $x\in \Tran_{k\mathcal{F}} (X, T)$. By Proposition~\ref{orbtranspoint}  $\orb_T (x)\subset \Tran_{k\mathcal{F}} (X, T)$, and hence $\Tran_{k\mathcal{F}} (X, T)$ is a dense $G_\delta$ subset of $X$ since $x\in \Tran (X, T)$. This finishes the proof.
\end{proof}

\begin{rem} \label{201604122105}
 Observe that when the state space $X$ is a compact metric space without isolated
points, $x\in \Tran (X, T)$ if and only if $x\in \Tran_\mathcal{B} (X, T)$. The family $\mathcal{F}_{\text{cof}}$ is clearly translation invariant (and hence $+$ invariant) and admits a countable base, and $k \mathcal{F}_{\text{cof}}= \mathcal{B}$. Thus by Proposition \ref{prop0201} one has: $(X, T)$ is transitive if and only if $\Tran (X, T)$ is a dense
    $G_\delta$ subset of $X$ if and only if $\Tran (X, T)\neq \varnothing$.
\end{rem}

\section{Transitive sensitivity and sensitive compactness}

 Recall that a dynamical
system $(X, T)$ is transitively sensitive if there exists  $\delta> 0$ such that  $S_T (W, \delta)\cap N_T (U, V)\neq \varnothing$ for
any opene subsets $U, V, W$ of $X$; and
sensitive compact if there exists $\delta> 0$ such that for any point $x\in X$ the  set
$\omega_{\mathcal{S}_T(\delta)}(x)$ is nonempty. Sometimes in that cases we will say also $(X, T)$ is transitively sensitive with a sensitive constant
$\delta$ and $(X, T)$ is sensitive compact with a sensitive constant
$\delta$. The main result of this section is the following

\begin{thm} \label{201604092203}
Let $(X, T)$ be a minimal system.  Then
the following conditions are equivalent:
\begin{enumerate}

\item $(X, T)$ is multi-sensitive.

\item $(X, T)$ is sensitive compact.

\item There exists $\delta > 0$ such that $\omega_{\mathcal{S}_T (\delta)} (x)= X$ for each $x\in X$.

\item There exist $\delta > 0$ and $x\in X$ with $\omega_{\mathcal{S}_T (\delta)} (x)= X$.

\item $(X, T)$ is transitively sensitive.
\end{enumerate}
 \end{thm}

 Before proceeding, we need:

\begin{lem} \label{201604092124}
Let $\delta> 0$ and $x\in X$. If $T: X\rightarrow X$ is almost open, then the family $\mathcal{S}_T (\delta)$ is $-$ invariant, and the subset $\omega_{\mathcal{S}_T (\delta)} (x)$ is positively $T$-invariant.
\end{lem}

\begin{proof}
By Proposition \ref{201604161045} it suffices to prove that $\mathcal{S}_T (\delta)$ is a $-$ invariant family.
Take arbitrary $F\in \mathcal{S}_T (\delta)$ and any $i\in \mathbb{Z}_+$. Then there exists opene subset $U$ of $X$ with $S_T (U, \delta)\subset F$. As $T: X\rightarrow X$ is almost open, $T^i: X\rightarrow X$ is also almost open, and then we can choose opene $V \subset T^i U$.
One has $g^{- i} (F)\supset g^{- i} S_T (U, \delta)= S_T (T^i U, \delta)\supset S_T (V, \delta)$, which implies that the family $\mathcal{S}_T (\delta)$ is $-$ invariant.
\end{proof}

The following result gives a characterization of transitive sensitivity for a general dynamical system in terms of dynamical compactness.

\begin{prop} \label{201604092152}
Let $(X, T)$ be a dynamical system.  Then the family $\mathcal{S}_T (\delta)$ is $+$ invariant for any $ \delta> 0$.
Furthermore, the following
conditions are equivalent:
\begin{enumerate}

\item $(X, T)$ is transitively sensitive.

\item There exist a $\delta > 0$ and  a dense $G_\delta$ subset $X_0\subset X$ such that $\omega_{\mathcal{S}_T (\delta)} (x)= X$ for each $x\in X_0$.

\item There exist a $\delta > 0$ and  a point $x\in X$ with $\omega_{\mathcal{S}_T (\delta)} (x)= X$.

\end{enumerate}
\end{prop}

\begin{proof}
Firstly, we show that $\mathcal{S}_T (\delta)$ is a $+$ invariant family. In fact, take any $F\in \mathcal{S}_T (\delta)$ and
each $i\in \mathbb{Z}_+$. We choose opene subsets $U, V$ of $X$ with $F\supset S_T (U, \delta)$ and $V\subset T^{- i} U$
satisfying $\diam (T^j V)< \delta$ for all $j= 0, 1, \dots, i$. Thus $g^i (F)\supset g^i S_T (U, \delta)\supset S_T (V, \delta)$
from the construction, and then $g^i (F)\in \mathcal{S}_T (\delta)$. This implies the $+$ invariance of the family $\mathcal{S}_T (\delta)$.

Observe that $(X, T)$ is transitively sensitive with a sensitive constant $\delta$, if and only if $(X, T)$ is
$k \mathcal{S}_T (\delta)$-transitive; and
that the family $\mathcal{S}_T (\delta)$ has a countable base: let $\mathcal{U}$ be a countable base of the family
of all opene subsets of $X$, then $\{S_T (U, \delta): U \in \mathcal{U}\}$ is a countable base
of $\mathcal{S}_T (\delta)$.
Then applying Proposition \ref{prop0201} the equivalence of $(1)\Leftrightarrow (2)\Leftrightarrow (3)$ follows from the fact that $x\in \Tran_{k \mathcal{S}_T (\delta)} (X, T)$ if and only if
$\omega_{\mathcal{S}_T (\delta)} (x)= X$.
\end{proof}

Observe that by \cite[Corollary 1.7]{GlasnerWeiss1993} the sensitivity of a dynamical system can be lifted up from a factor to an extension by an almost open factor map between transitive systems. The following result gives the lift-up and projection property of transitive sensitivity between transitive systems.

\begin{lem} \label{201604101854}
Let $\pi: (X, T)\rightarrow (Z, R)$ be a factor map between dynamical systems.
\begin{enumerate}

\item Assume that $\pi$ is almost one-to-one. If $(Z, R)$ is transitively sensitive with a sensitive constant $\delta> 0$ then $(X, T)$ is also transitively sensitive.

\item Assume that there exists $z\in Z$ whose fiber is a singleton. If $(X, T)$ is transitively sensitive then $(Z, R)$ is sensitive, in particular, $\text{Eq} (Z, R)= \varnothing$.
\end{enumerate}
\end{lem}

\begin{proof}
(1) We take a compatible metric $\rho$ over $Z$ and let $\varepsilon> 0$ such that $d (x_1, x_2)\le \varepsilon$
implies $\rho (\pi x_1, \pi x_2)\le \delta$ for any $x_1, x_2\in X$. Now let $U, V, W$ be arbitrary opene subsets of $X$.
As the map $\pi: X\rightarrow Z$ is almost one-to-one, we may take opene subsets $U_Z, V_Z, W_Z$ of $Z$
with $\pi^{- 1} U_Z\subset U$, $\pi^{- 1} V_Z\subset V$ and $\pi^{- 1} W_Z\subset W$. Observe that
if $n\in S_R (W_Z, \delta)\cap N_R (U_Z, V_Z)$, then: on one hand, there exist $z_1, z_2\in W_Z$ with $\rho (R^n z_1, R^n z_2)> \delta$,
and so $d (T^n x_1, T^n x_2)> \varepsilon$ for any $x_1\in \pi^{- 1} (z_1)$ and $x_2\in \pi^{- 1} (z_2)$, hence $\diam (T^n W)> \varepsilon$;
on the other hand, $U_Z\cap R^{- n} V_Z\neq \varnothing$, and then $U\cap S^{- n} V\supset \pi^{- 1} U_Z\cap \pi^{- 1} (R^{- n} V_Z)\neq \varnothing$. This implies $S_T (W, \varepsilon)\cap N_T (U, V)\supset S_R (W_Z, \delta)\cap N_R (U_Z, V_Z)\neq \varnothing$, as $(Z, R)$ is transitively sensitive. Thus, by the
arbitrariness of $U, V$ and $W$, we have that $(X, T)$ is also transitively sensitive.

(2) As $(X, T)$ is transitively sensitive (and assume with a sensitive constant $\delta> 0$), it is clear that $(Z, R)$ is transitive, and
then
by the refined Auslander-York dichotomy the system $(Z, R)$ is sensitive if and only if $\text{Eq} (Z, R)= \varnothing$
(see \cite{AuslanderYorke}, \cite{GlasnerWeiss1993}, \cite{AAB1993} and the book \cite{Akin1997}). Thus it suffices to prove
$\text{Eq} (Z, R)= \varnothing$. Let $\rho$ be a compatible metric over $Z$, and assume the contrary to take a point $z\in \text{Eq} (Z, R)$.
By the assumption that there exists a point of $Z$ whose fiber is a singleton, we may take an opene subset $W$ in $Z$ with
$\diam (\pi^{- 1} W)< \delta$, and an opene subset $W_*\subset W$ and $\delta_*> 0$ such that if the distance between a point of $Z$
and $W_*$ is smaller than $\delta_*$ then the point belongs to $W$. Since $z\in \text{Eq} (Z, R)$, there exists an open neighborhood
$U_*$ of $z$ with $\diam (R^n U_*)< \delta_*$ for all $n\in \mathbb{Z}_+$. As $(X, T)$ is transitively sensitive with a sensitive
constant $\delta$, take $m\in N_T (\pi^{- 1} U_*, \pi^{- 1} W_*)\cap S_T (\pi^{- 1} U_*, \delta)$. Thus $T^m (\pi^{- 1} U_
*)\cap \pi^{- 1} W_*\neq \varnothing$, and then $R^m U_*\cap W^*\neq \varnothing$, which implies $R^m U_*\subset W$ by the construction of
$U_*, W_*$ and $W$. Hence
$$\diam (T^m (\pi^{- 1} U_*))\le \diam (\pi^{- 1} (R^m U_*))\le \diam (\pi^{- 1} W)< \delta,$$
a contradiction to the selection of $m\in S_T (\pi^{- 1} U_*, \delta)$. This finishes the proof.
\end{proof}

 Now we are ready to prove Theorem \ref{201604092203}.

\begin{proof}[Proof of Theorem \ref{201604092203}]
$(1)\Rightarrow (2)$ follows directly from the definitions.
As the system $(X, T)$ is minimal, the map $T: X\rightarrow X$ is almost open. Observing that $\omega_{\mathcal{S}_T (\delta)} (x)$ is a closed subset of $X$ for each $x\in X$, the implication of
$(2)\Rightarrow (3)$ follows from Lemma \ref{201604092124} and the minimality of $(X, T)$. The implication of $(3)\Rightarrow (4)\Rightarrow (5)$ follows from Proposition \ref{201604092152}.
Since a minimal system is either multi-sensitive or a weakly almost one-to-one extension of its maximal equicontinuous factor by \cite{HKZ}, and then
$(5)\Rightarrow (1)$ follows from Lemma \ref{201604101854}. This finishes the proof.
\end{proof}

 Clearly each multi-sensitive system is sensitive compact. Observe that
each non-proximal, transitive compact system is multi-sensitive by \cite[Theorem 4.7]{HKhKZh2015}.
In particular, each minimal transitive compact system is multi-sensitive,
as each minimal proximal system is trivial by \cite{AK} and all dynamical systems considered are assumed to be nontrivial. Nevertheless, there are many minimal, non transitive compact, multi-sensitive systems.
For example, consider the classical dynamical system $(X, T)$ given by $X= \mathbb{R}^2/\mathbb{Z}^2$ and $T: (x, y)\mapsto (x+ \alpha, x+ y)$ with $\alpha\notin \mathbb{Q}$ (see \cite[Chapter 1]{Furstenberg1981}). As commented in \cite[Page 1816]{HKhKZh2015}, $(X, T)$ is an invertible minimal multi-sensitive system; note that $(X, T)$ is not weakly mixing, since $(X, T)$ admits an irrational rotation as its nontrivial equicontinuous factor and any equicontinuous factor of a weakly mixing system is trivial.
Remark that by \cite[Corollary 3.10]{HKhKZh2015} for a minimal system the system is transitive compact if and only if it is weakly mixing, and then the constructed system $(X, T)$ is not transitive compact.

\begin{prop} \label{201604181936}
Each weakly mixing system $(X, T)$ is transitively sensitive.
\end{prop}

\begin{proof}
Observe that we are only interested a nontrivial dynamical system, and then let $0< \delta< \diam (X)$. We choose opene subsets $W_1, W_2$ of $X$ such that the distance between $W_1$ and $W_2$ is strictly larger than $\delta$. Now take arbitrary opene subsets $U, V, W$ of $X$.
As $(X, T)$ is weakly mixing, $(X^3, T^{(3)})$ is transitive by \cite{Furstenberg1967}, and then $N_T (U, V)\cap S_T (W, \delta)\supset N_{T^{(3)}} (U\times W\times W, V\times W_1\times W_2)\neq \varnothing$. This implies that $(X, T)$ is transitively sensitive with a sensitive constant $\delta> 0$.
\end{proof}

We give a sufficient condition for a dynamical system being transitively sensitive (by Proposition \ref{201604092152}) as the end of this section.

\begin{lem}
Assume $\omega_{\mathcal{S}_T (\varepsilon)} (x)= X$ for some $x \in X$ and $\varepsilon>0$. Then there is $\delta> 0$ such that for any opene subset $U$ of $X$ and each
neighbourhood $U_x$ of $x$ there are $y\in U_x$ and $n\in N_T (x, U)$ with $d(T^n x, T^n y)>\delta$.
If in addition, the map $T: X\rightarrow X$ is almost one-to-one, then the converse holds.
\end{lem}

\begin{proof}
Fix an opene subset $U$ of $X$ and a neighborhood $U_x$ of $x$. As $\omega_{\mathcal{S}_T (\varepsilon)} (x)= X$, there is $n\in N_T (x, U)\cap S_T (U_x, \varepsilon)$, and then there are points
$x_1, x_2 \in U_x$ with $d(T^n x_1, T^n x_2)> \varepsilon$. We have either
$d (T^n x, T^n x_1)> \frac{\varepsilon}{2}$ or $d (T^n x, T^n x_2)> \frac{\varepsilon}{2}$, and then obtain the desired statement
for $\delta= \frac{\varepsilon}{2}$.

Now suppose
that there is $\delta> 0$ such that for any opene subset $U$ of $X$ and each
neighbourhood $U_x$ of $x$ there are $y\in U_x$ and $n\in N_T (x, U)$ with $d(T^n x, T^n y)>\delta$, and that the map $T: X\rightarrow X$ is almost one-to-one. Let $U, W$ be arbitrary opene subsets of $X$. It is clear $x\in \Tran (X, T)$, and so there is $k\in N_T (x, W)$. Note that $T: X\rightarrow X$ is almost one-to-one, the map $T^k: X\rightarrow X$ is also almost one-to-one by Corollary \ref{201604102314}, and then we may take an opene subset $V$ of $X$ with $T^{- k} V= (T^k)^{- 1} V\subset U$. Observe that $T^{- k} W$ is a neighborhood of $x$, and then by the assumption there exists a point $y\in T^{- k} W$ and an integer $n\in \mathbb{Z}_+$ such that $n\in N_T (x, V)$ and $d (T^n x, T^n y)> \delta$. In fact, we may assume $n>k$, else
we can replace $T^{- k} W$ by a small enough open neighbourhood $G_x\subset T^{- k} W$ of $x$ with $\diam (T^i G_x)< \delta$ for all $0\le i\le k$.
Then $\diam (T^{n- k} W)> \delta$ as $T^k x\in W$ and $T^k y\in W$, and $T^{n- k} x\in T^{- k} (T^n x)\subset T^{- k} V\subset U$. In particular, $N_T (x, U)\cap S_T (W, \delta)\neq \varnothing$.
Thus $\omega_{\mathcal{S}_T (\delta)} (x)= X$ by the arbitrariness of $U$ and $W$.
\end{proof}

 \section{Weakly disjointness and weakly mixing}

Recall that dynamical systems $(X, T)$ and $(Y, S)$ are \emph{weakly disjoint} if the product system $(X\times Y, T\times S)$ is transitive. The following theorem characterizes weak disjointness,
which is proved firstly by Weiss \cite{Weiss2000} in some special class and
then is generalized by Akin and Glasner \cite{AkinGlasner2001}.
We say that $\mathcal{F}$ is \emph{thick} if $\tau \mathcal{F}= \mathcal{F}$, where
$$\tau \mathcal{F}= \left\{F\subset \mathbb{Z}_+: \bigcap_{j= 1}^n g^{- i_j} F\in \mathcal{F}\ \text{for each}\ n\in \mathbb{N}\ \text{and all}\ i_1, \dots, i_n\in \mathbb{Z}_+\right\}.$$

\begin{thm*}
Let $\mathcal{F}$ be a proper, translation invariant, thick family. A dynamical system is
$k\mathcal{F}$-transitive if and only if it is weakly disjoint from every $\mathcal{F}$-transitive system.
\end{thm*}

Observe that a dynamical system is
weakly mixing if and only if it is weakly disjoint from itself, and then weak disjointness is
characterized by \cite[Proposition 3.8]{HKhKZh2015} in some special case. Now we discuss weak disjointness using
dynamical compactness which will be some generalization of  \cite[Proposition 3.8]{HKhKZh2015}. We will need the following

\begin{lem} \label{lem0302}
 Let $(X, T)$ and $(Y, S)$ be dynamical systems and let $x\in X$.
Then the family $\mathcal{N}_S$ is translation invariant and $\omega_{\mathcal{N}_S} (x)= \omega_{\mathcal{N}_S} (T x)$.
\end{lem}

\begin{proof}
By Proposition \ref{201604161045} it suffices to prove that $\mathcal{N}_S$ is a translation invariant family. We also suppose that $\mathcal{N}_S$ is proper (i.e., $(Y,S)$ is a transitive system) since, otherwise, the result is trivial.
Take arbitrary $F\in \mathcal{N}_S$ and any $i\in \mathbb{Z}_+$. Then there exist opene subsets $U, V$ of $Y$ with $N_S (U, V)\subset F$.
As the non-singleton space $Y$ contains no isolated points, we can take  suitable opene
$V_1\subset V$ and $U_1\subset U$ such that
$U_1\cap \bigcup_{k= 0}^i S^k \overline{V_1}=\varnothing $.
 One has $g^{i} N_S (U, V)\supset N_S (S^{- i} U_1, V_1)$, which implies that the family $\mathcal{N}_S$ is $+$ invariant: in fact, if $n\in N_S (S^{- i} U_1, V_1)$ then $n> i$ by the selection, and so $n- i\in N_S (U_1, V_1)\subset N_S (U, V)$. Moreover, it is clear  $g^{- i} (F)\supset g^{- i} N_S (U, V)\supset N_S (U, S^{- i} V)$, and then the family $\mathcal{N}_S$ is $-$ invariant. This finishes the proof.
\end{proof}

\begin{thm} \label{th0201}
The following conditions are equivalent:
\begin{enumerate}

\item The systems $(X, T)$ and $(Y, S)$ are weakly disjoint.

\item
Both $\Tran_{k \mathcal{N}_S} (X, T)$ and $\Tran_{k \mathcal{N}_T} (Y, S)$ are dense $G_\delta$ subsets.

\item
The set $\Tran_{k \mathcal{N}_S} (X, T)$ is a dense $G_\delta$ subset of $X$.

\item
Both $\Tran_{k \mathcal{N}_S} (X, T)$ and $\Tran_{k \mathcal{N}_T} (Y, S)$ are nonempty subsets.

\item
The set $\Tran_{k \mathcal{N}_S} (X, T)$ is a nonempty subset of $X$.
\end{enumerate}
\end{thm}

\begin{proof}
$(1)\Leftrightarrow (2)\Leftrightarrow (3)$:
It is clear from the definition that: the system $(X, T)$ is $k \mathcal{N}_S$-transitive, if and only if the systems $(X, T)$ and $(Y, S)$ are weakly disjoint, if and only if $(X, T)$ is $k \mathcal{N}_S$-transitive and $(Y, S)$ is $k \mathcal{N}_T$-transitive. As both $\mathcal{N}_T$ and $\mathcal{N}_S$ are families admitting a countable base, it is direct to obtain the equivalence of $(1)\Leftrightarrow (2)\Leftrightarrow (3)$ by applying Proposition \ref{prop0201}.

The implication $(2)\Rightarrow (4)\Rightarrow (5)$ is obvious. To finish the proof, we only need to show $(5)\Leftrightarrow (3)$. By Lemma \ref{lem0302} the family $\mathcal{N}_S$ is translation invariant, and then  the family $k \mathcal{N}_S$ is also translation invariant. Thus the equivalence of $(5)\Leftrightarrow (3)$ follows from Proposition \ref{prop0201}.
\end{proof}

Note that we have a characterization of weak mixing by using dynamical compactness \cite[Proposition 3.8]{HKhKZh2015}. Now we
improve \cite[Proposition 3.8]{HKhKZh2015} as follows.

  Recall that $\mathcal{S}\subset \mathbb{N}$ is an \emph{IP set}
 if there exists $\{p_k: k\in \mathbb{N}\}\subset \mathbb{N}$ with $FS\{p_i\}_{i=1}^\infty\subset \mathcal{S}$, where
$FS\{p_i\}_{i=1}^\infty=\{p_{i_1}+ \dots+ p_{i_k}: k\in \mathbb{N}\
 \text{and}\ 1\le i_1< \dots< i_k\}.$
Analogously, for each $n\in \mathbb{N}$ we define
$FS\{p_i\}_{i=1}^n= \{p_{i_1}+ \dots+ p_{i_k}: k\in \mathbb{N}\
 \text{and}\ 1\le i_1< \dots< i_k\le n\}.$
 Denote by $\mathcal{F}_{ip}$ the family of all IP sets.

By \cite[Theorem 3.2]{Li}, the subset $\Tran_{\mathcal{F}_{ip}} (X, T)$
contains a dense $G_\delta$ subset of $X$ for any weakly mixing system $(X, T)$,
while $\Tran_{\mathcal{F}_{ip}} (X, T) \neq \varnothing$ does not imply the weak mixing of the system $(X, T)$ by \cite[Proposition 3.4]{Li}.
We will improve that in the following Proposition \ref{0203}. Before proceeding, we make the following

\begin{lem} \label{l4.1}
Let $(X, T)$ be a dynamical system and $\mathcal{F}$ be a family.

\begin{enumerate}
\item
Let $\delta> 0$. If the family $\mathcal{S}_T (\delta)\cdot \mathcal{F}$ is proper then $\mathcal{S}_T (\delta)\cdot \mathcal{F}\subset \mathcal{B}$.

\item If the family $\mathcal{N}_T\cdot \mathcal{F}$ is proper then $\mathcal{N}_T\cdot \mathcal{F}\subset \mathcal{B}$.
\end{enumerate}
\end{lem}

\begin{proof}
(1) Assume the contrary that there exists an opene subset $U$ in $X$ and $F \in \mathcal{F}$ such that $S_T (U, \delta)\cap F$ is finite, and so we may choose $m \in \mathbb{N}$ such that $n \notin S_T(U, \delta)\cap F$ for any integer $n>m$. Since $T:X \rightarrow X$ is uniformly continuous one can find opene $V \subset U$ small enough such that $\diam (T^i V) < \delta$ for all $0\leq i \leq m$.
Then $S_T (V, \delta) \subset S_T (U, \delta)$, which implies $S_T(V, \delta) \cap F= \varnothing$. A contradiction.

(2)
Assume the contrary that there exist opene subsets $U, V$ in $X$ and $F \in \mathcal{F}$ such that $N_T(U, V)\cap F$ is finite,
say $N_T(U, V)\cap F= \{n_1, \dots, n_k\}$. As the non-singleton space $X$ contains no isolated points, we can take opene
 $U_1\subset U$ small enough such that
$V_1:= V\setminus \bigcup_{i= 1}^k T^{n_i} \overline{U_1}$ is an opene subset of $X$. By the construction we
have $N_T(U_1, V_1)\subset N_T (U, V)$ and then $N_T(U_1, V_1)\cap F= \varnothing$, a contradiction.
\end{proof}

\begin{prop} \label{0203}
The following conditions are
equivalent:
\begin{enumerate}
\item $(X, T)$ is weakly mixing.

\item There exists a dense $G_\delta$ subset $X'$ of
$X$ such that, for each $x\in X'$, $   n_T(x, G) \cap N_T(U, V) \in \mathcal{F}_{ip}$ for all opene subsets $G, U, V$ of $X$.
\end{enumerate}
\end{prop}

\begin{proof}
$(2)\Rightarrow (1)$: Just observe from the assumption that $\omega_{\mathcal{N}_T}(x)=X$ for all $x \in X'$, and hence the system $(X, T)$ is weakly mixing by \cite[Proposition 3.8]{HKhKZh2015}.

(1) $\Rightarrow$ (2): Since $(X, T)$ is weakly mixing, $(X^2, T^{(2)})$ is also weakly mixing by \cite{Furstenberg1967} and \cite{Petersen}, and hence by \cite[Proposition 3.8]{HKhKZh2015} there is a dense $G_\delta$ subset
$Y \subset X^2$ such that
$\omega_{\mathcal{N}_{T^{(2)}}}((x_1, x_2))=X^2$ for each $(x_1, x_2) \in Y$. Applying the well-known Ulam Lemma
there is a dense $G_\delta$ subset $X' \subset X$ such that, for any $x \in X'$, $\{y: (x, y) \in Y\}$ is a dense $G_\delta$ subset of $X$. Now we show that $X'$ is the
desired set.

Let $x \in X'$ and fix any opene subsets $G, U, V$ of $X$. Choose $y \in G$ with $(x, y) \in Y$ and then $\omega_{\mathcal{N}_{T^{(2)}}}((x, y))=X^2$, in particular,
$(y, y) \in \omega_{\mathcal{N}_{T^{(2)}}}((x, y))$. Thus
 $$n_T(x, G)\cap  n_T(y, G) \cap N_T(U, V) \cap N_T(V, V) \neq \varnothing,$$
and take $p_1\in \mathbb{N}$ from this set by Lemma \ref{l4.1}. We have $p_1 \in    n_T(x, G) \cap N_T(U, V)$ and
$T^{p_1} y \in G$, $T^{p_1} V \cap V \neq \varnothing$. Define opene subsets $G_1= G \cap T^{- p_1}G\ni y$ and $V_1=V \cap T^{-p_1}V$.
Now we proceed inductively. Suppose that we are given a
sequence $\{p_1, \dots, p_k\}\subset \mathbb{N}$ with $FS\{p_i\}_{i=1}^{k} \subset    n_T(x, G) \cap N_T(U, V)$, and opene subsets
$$G_k=\bigcap_{s \in FS\{p_i\}_{i=1}^k \cup \{0\}} T^{-s} G\ni y, V_k=\bigcap_{s \in FS\{p_i\}_{i=1}^k \cup \{0\}} T^{-s}V.$$
As $(y, y) \in \omega_{\mathcal{N}_{T^{(2)}}}((x, y))$, we may take $p_{k+1}\in \mathbb{N}$ by Lemma \ref{l4.1} from the set
$n_T(x, G_k) \cap  n_T(y, G_k) \cap N_T(U, V_k) \cap N_T(V_k, V_k)$. It is not hard to check that
 $$G_{k+1}=G_k \cap T^{-p_{k+1}}G_k=\bigcap_{s \in FS\{p_i\}_{i=1}^{k+ 1} \cup \{0\}} T^{-s} G\ni y,$$
 $$V_{k+1}=V_k \cap T^{-p_{k+1}}V_k=\bigcap_{s \in FS\{p_i\}_{i=1}^{k+1} \cup \{0\}} T^{-s}V$$
 are both opene subsets of $X$, and that
 $FS \{p_i\}_{i=1}^{k+1} \subset    n_T(x, G) \cap N_T(U, V)$, which completes the induction.
Finally, $n_T(x, G) \cap N_T(U, V)\in \mathcal{F}_{ip}$ with $FS \{p_i\}_{i=1}^{\infty} \subset n_T(x, G) \cap N_T(U, V)$. This finishes the proof.
\end{proof}

\section{Transitive compact (non weakly mixing) systems}

     Recall that the system $(X, T)$ is \emph{totally transitive} if $(X, T^k)$ is transitive for each $k\in \mathbb{N}$; and is \emph{topologically mixing} if $N_T(U,V)\in \mathcal{F}_\text{cof}$ for any opene subsets $U,V$ in $X$.
Note that $(X, T)$ is weakly mixing if and only if
$N_T(U,V)\in \mathcal{F}_{\text{thick}}$ for any opene sets $U,V$ in $X$ by \cite{Furstenberg1967, Petersen}, and so any weakly mixing system is totally transitive.
It is direct to check that each weakly mixing system is transitive compact. In \cite{HKhKZh2015} the authors showed the existence
of non totally transitive, transitive compact systems in both proximal and non-proximal cases. We extend it as follows:

\begin{thm} \label{th0601}
There is a totally transitive, transitive compact system $(X, T)$ which is not weakly mixing.
\end{thm}

\begin{proof}
  Take a nontrivial proximal, topologically mixing system $(Y, S)$ and let $(\mathbb{S}^1, R_\alpha)$ be the standard irrational rotation
 on the unit circle $\mathbb{S}^1=\mathbb{R}/\mathbb{Z}$ with $\alpha\notin \mathbb{Q}$. Note that a dynamical system is
proximal if and only if it contains the unique fixed point, which is the only minimal point of the system \cite{AK}.
Denote by $p_Y$ the unique minimal point (fixed point) of $(Y, S)$. Observe that the system $(Y \times \mathbb{S}^1, S\times R_\alpha)$ is
totally transitive: for each $n\in \mathbb{N}$, the system $(Y, S^n)$ is topologically mixing by the definition and it is standard that the system $(\mathbb{S}^1, R_\alpha^n)$ is minimal, then it is direct to see that these two systems are weakly disjoint.

 Let $(X,T)$ be the quotient system $Y \times \mathbb{S}^1/ \sim$ equipped with the action $T$ induced naturally from
 $S\times R_\alpha$, where the equivalence relation $\sim$ is defined via: given $x, y\in X$, $x\sim y$ if and only if
 either $x=y$ or $x$ and $y$ both have $p_Y$ in the first coordinate.
 In other words the space $X$ looks like a cone space,  where the vertex of the cone is a point $p$, each ``horizontal'' fiber spaces are the space $Y$, the vertical fiber spaces are
 the circles (see Figure 1). Clearly, $(X, T)$ is totally transitive.

 \medskip

\begin{center}
\usetikzlibrary{calc}

\begin{tikzpicture} [rotate=-90]
  \def\rx{2}    
  \def\ry{0.5}  
  \def\z{3}     

  \pgfmathparse{asin(\ry/\z)}
  \let\angle\pgfmathresult

  \coordinate (h) at (0, \z);
  \coordinate (O) at (0, 0);
  \coordinate (A) at ({-\rx*cos(\angle)}, {\z-\ry*sin(\angle)});
  \coordinate (B) at ({\rx*cos(\angle)}, {\z-\ry*sin(\angle)});

  \draw[fill=gray!50] (A) -- (O) -- (B) -- cycle;
  \draw[fill=gray!30] (h) ellipse ({\rx} and {\ry});
  \draw (0,-0.2) node{\footnotesize ~p};
  \draw (1.15, 1.25) node{Y};
   \draw (0,2.75) node{$\mathbb{S}^1$};
\end{tikzpicture}

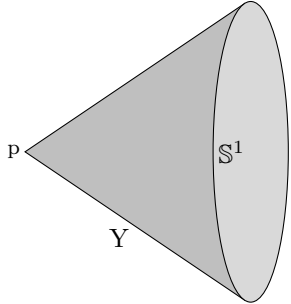
{\footnotesize Figure 1.}
\end{center}
\smallskip

Denote by $q:Y \times \mathbb{S}^1 \rightarrow X$ the corresponding quotient map, then $q: Y_\infty\times \mathbb{S}^1\rightarrow X \setminus \{p\}$ is a homeomorphism, where we set $Y_\infty= Y\setminus \{p_Y\}$.
It is standard that the system $(\mathbb{S}^1, R_\alpha)$ is not weakly mixing, and then there exist opene subsets $U_*, V_*$ of $\mathbb{S}^1$ with $N_{R_\alpha} (U_*, V_*)\notin \mathcal{F}_{\text{thick}}$, hence
$$N_{T} (q (Y_\infty\times U_*), q (Y_\infty\times V_*))= N_{S\times R_\alpha} (Y_\infty\times U_*, Y_\infty\times V_*)\subset N_{R_\alpha} (U_*, V_*)$$
is not thick. This implies that the system $(X, T)$ is not weakly mixing.

Now let $U, V$ be arbitrary opene subsets of $X$. We can choose opene subsets $U_1, V_1\subset Y_\infty$ and $U_2, V_2 \subset \mathbb{S}^1$
with $U_1\times U_2 \subset q^{- 1} U$ and $V_1\times V_2 \subset q^{- 1} V$. As $(Y, S)$ is topologically mixing and $(\mathbb{S}^1, R_\alpha)$ is minimal, $N_S (U_1, V_1)\in \mathcal{F}_\text{cof}$ and $N_{R_\alpha} (U_2, V_2)\in \mathcal{F}_\text{syn}$, and then $N_S (U_1, V_1)\cap N_{R_\alpha} (U_2, V_2)\in \mathcal{F}_\text{syn}$, thus
$$N_{T} (U, V)= N_{S\times R_\alpha} (q^{- 1} U, q^{- 1} V)\supset N_{S\times R_\alpha} (U_1\times U_2, V_1\times V_2)$$
is a syndetic set. Observe from the construction that the system $(X, T)$ is proximal with $p$ as its unique fixed point, then $n_T (x, U_p)$ is a thickly syndetic subset for each point $x\in X$ and any neighbourhood $U_p$ of $p$ (see \cite[Lemma 3.12]{HKhKZh2015}). This implies $p \in \omega_{\mathcal{N}_T} (x)$ for each $x\in X$, and then the system $(X, T)$ is transitive compact.
\end{proof}

The following result is proved independently in \cite{DoFr} and \cite{Sh}.

\begin{lem} \label{Shar}
Any $\omega$-limit set $\omega_T(x)$ can not be decomposed into $\alpha$ disjoint closed, nonempty, positively $T$-invariant subsets, where $2\le \alpha\le \aleph_0$.
\end{lem}

Before proceeding, we need the following example, for which we fail to find a reference and hence provide a detailed construction, as it is crucial in our arguments.

\begin{prop} \label{exam}
For any given compact metric space $Z$, there exists a topologically mixing system $(X, T)$ such that, $Z$ can be realized as the set of all its minimal points, furthermore, its each minimal point is a fixed point.
\end{prop}

\begin{proof}
The construction is divided into two steps.

\medskip

In the first step we shall construct a
 topologically mixing system $(Y, F)$ but with two fixed points, which are the only minimal points of the system.
 Let $\Sigma=\{0,1\}^{\mathbb{Z}_+}$ and
 $\sigma :\Sigma \to \Sigma$ be the full (one-sided) shift.
We are going to find the system $(Y, F)$ of the form $(\overline{\orb_{\sigma} (x)}, \sigma)$  for some $x \in \Sigma$.

In order to define $x\in \Sigma$, firstly we represent each $W\in \Sigma^*$ with $|W|\ge 1$ in the following form:
$W=a^{i} Q b^{j}$, where $a^{i}$ and $b^{j}$
 (with $i\ge 1$) are the longest segments of equal digits which we can take at the beginning and at the end of $W$,
whereas $Q$ is the rest, possible the empty subblock. Clearly, $j$ may be equal to $0$ and then $b^j$ will be the empty subblock, in this case, we treat the digit $b$ as $0$; in particular, if $W= a^k$ then we set $Q$ to be the empty subblock and $i= k$, $b= 0$, $j= 0$.

Now we are going to define $x\in \Sigma$. Let $A_1=10$ be the first block of $x$ and define inductively
the rest blocks $A_2, A_3, \dots$,  then $x$ will be the limit of the starting blocks $A_k$.
Suppose that we have defined $A_k, k\in \mathbb{N}$. Since $A_{k}$ has finitely many subblocks, there is a finite number of different
pairs of these subblocks.
For any pair $(W_1, W_2)$ of subblocks of $A_k$ we will define a block $c (W_1, W_2)$ by using their combination.
Then we are ready to define $A_{k+1}$: at the beginning of $A_{k+1}$ we write $A_k0^k1^k$, and then all possible blocks $c (W_1, W_2)$ of
pairs $(W_1, W_2)$ of subbloks of $A_{k}$ in any fixed order.
The definition $c (W_1, W_2)$ depends on the structure of $W_1$ and $W_2$. Let us write $W_1$ and $W_2$ in the form as above:
$W_1=a^{i_1} Q_1 b^{j_1}$ and $W_2=c^{i_2} Q_2 d^{j_2}$,
 where $a^{i_1}$, $b^{j_1}$, $c^{i_2}$ and $d^{j_2}$ (with $i_1, i_2\ge 1$)
are the longest segments of equal digits which we can take at the beginning and at the end of $W_1$ and $W_2$,
whereas $Q_1$ and $Q_2$ are the rest, possible empty subblocks ($j_1$ and $j_2$ may be equal to $0$, and then we treat the
corresponding digits as $0$). \emph{The combination block} of the pair $W_1, W_2$, i.e. $c(W_1, W_2)$, is defined as follows:
$$c(W_1, W_2)=a^{k+i_1} Q_1 b^{j_1+k}c^{k+i_2}Q_2d^{k+j_2}a^{k+i_1} Q_1 b^{j_1+k+1}c^{k+i_2}Q_2 d^{k+j_2} .$$
We see from the construction that $A_{k+1}$ is presented as a sequence of blocks with length not longer
than $|A_k|$, which are separated from each other with some
sequences of blocks of consecutive 1's or 0's of length not less than $k$.
In fact, for all
$m>k$ this property holds for $A_m$ (and hence in any subblock of $A_m$ with length more than $|A_k|+ 2k- 1$ one can find $0^k$ or $1^k$).
Suppose that $A_m$ may be presented in this form.
Observe that
$A_{m+1}$ is obtained by adding to $A_m 0^m 1^m$
combination blocks $c(W_1^*, W_2^*)$ of all pairs $(W_1^*, W_2^*)$ of subblocks of $A_m$ glued in a proper way, and note that if the property holds for blocks $W_1^*$ and $W_2^*$ (with beginning and end parts of the blocks possible exception) then the property holds
for the combination block $c(W_1^*, W_2^*)$. In particular, the property holds for $A_{m+ 1}$.

Put $(Y, F):=(\overline{\orb_{\sigma} (x)},\sigma)$. Now let us check that it has the required properties.

Firstly, in any subblock of $x$ with length more than $|A_k|+2k-1$ one can find $0^k$ or $1^k$. And hence, if we increase the length of a block $w$ of $x$ the length of the biggest subblock of $w$ in form of $0^m$ or $1^m$ increases unboundedly. Moreover, for each $k\in \mathbb{N}$ both $0^k$ and $1^k$ appear in $x$.
In particular,
both
$0^\infty$ and $1^\infty$ are fixed points of $(\overline{\orb_{\sigma} (x)}, \sigma)$, and there is no other minimal sets in it.

Recall that the base for the open sets in $\Sigma$ is given by the collection of all cylinder sets
$C[c_0c_1c_2 \dots c_m]=\{x \in \Sigma:~ x_i=c_i \text{ for } 0\le i\leq m \}.$
 Given any two subblocks $W_1$ and $W_2$ of $A_m$, we write $W_1$ and $W_2$ in the form as above:
$W_1=a^{i_1} Q_1 b^{j_1}$ and $W_2=c^{i_2} Q_2 d^{j_2}$. If $b= c$, then for each $l\ge m$ there exists a combination block $c_l (W_1, W_2)$
 containing the subblock $W_1 b^l W_2$, and hence $N_F ([W_1]\cap Y, [W_2]\cap Y)\supset \{m+ |W_1|, m+ |W_1|+ 1, \dots\}$;
if $b\neq c$, then for each $l\ge m$ there exists a combination block $c_l (W_1, W_2)$
 containing the subblocks $W_1 b^k c^k W_2$ and $W_1 b^k 1 c^k W_2$, and hence $N_F ([W_1]\cap Y, [W_2]\cap Y)\supset \{2 m+ |W_1|, 2 m+ |W_1|+ 1, \dots\}$. This shows that the system $(Y, F)$ is topologically mixing.

\medskip

Now we shall finish the construction by the second step. Firstly, we take $(X', T')$ to be the product system $\prod_1^\infty (Y, F)$. It is ready to check that the system $(X', T')$ is topologically mixing, for which the middle-third Cantor set $C$ is the set of all its minimal points and its each minimal point is a fixed point. Note that there exists a continuous surjection $h: C\rightarrow Z$ (see for example \cite[Page 165-166, Problem O]{Kelley}), and then we consider the quotient system $(X, T)$ with $X= X'/\sim$ equipped with the action induced naturally from $T'$, where the closed positively $T'\times T'$-invariant equivalence relation $\sim$ is defined via $x\sim y$ if and only if $x= y\in X'\setminus Z$ or $h (x)= h (y)$ for $x, y\in Z$. Then the system $(X, T)$ has the required properties.
\end{proof}

The following result shows that in general there is no a topological structure similar to Lemma \ref{Shar} for the $\omega_{\mathcal{N}_T}$-limit sets.

\begin{thm} \label{20151026}
For any given compact metric space $Z$, there exists a non totally transitive, transitive compact system $(X, T)$ such that,
 $Z$ can be realized as the set of all its minimal points with its each minimal point being a fixed point, furthermore,
 $Z$ is realized as $\omega_{\mathcal{N}_T} (x)$ for some $x\in X$.
\end{thm}

\begin{proof}
 The idea of the proof is very similar to that of the first part (proximal case) of \cite[Theorem 3.14]{HKhKZh2015}.
 Instead of
 a nontrivial proximal, topologically mixing system $(Y, F)$ there
 (main point is ``a map with exactly one minimal point!''), we take again a
 topologically mixing system $(Y, F)$, but with $Z$ realized as the set of all its minimal points where its each minimal point is a fixed point (for existence of such a dynamical system see Proposition \ref{exam}). Then, by the wedge sum construction there, we obtain a non totally transitive system $(X, T)$ such that $Z$ is realized as the set of all its minimal points with its each minimal point being a fixed point. Similar to arguments there, it is not hard to show $Z\supset \omega_{\mathcal{N}_T} (x)$ for all $x\in X$.

Firstly we prove the following claim:

 \medskip

 \noindent {\bf Claim.} For $x\in X$, if $x= y\in Y$ and $p\in Z\cap \omega_{\mathcal{N}_{F^2}} (y)$ then, as a point in $X$, $p$ belongs to $\omega_{\mathcal{N}_T} (x, T)$.

 \begin{proof}[Proof of Claim]
Let $U_p$ be an open subset of $X$ containing $p$, and clearly $U_p$ may be also viewed as an open subset of $Y$ containing $p$. Now for any given opene subsets $U$ and $V$ of $X$: if both $U$ and $V$ can be viewed as opene subsets of $Y$, then we can take $n\in N_{F^2} (y, U_p)\cap N_{F^2} (U, V)$ by the assumption $p\in \omega_{\mathcal{N}_{F^2}} (y)$ and hence $2 n\in N_{T} (y, U_p)\cap N_{T} (U, V)$; if both $U$ and $V$ can be viewed as opene subsets of $Y_c$, then both $T^{- 1} U$ and $T^{- 1} V$ can be viewed as opene subsets of $Y$ and hence $N_{T} (y, U_p)\cap N_{T} (U, V)= N_{T} (y, U_p)\cap N_{T} (T^{- 1} U, T^{- 1} V)\neq \varnothing$; if $U$ and $V$ can be viewed as opene subsets of $Y$ and $Y_c$, respectively, noting $p\in Z$ and hence $T p= p$, there is an opene subset $V_p$ of $X$ containing $p$ such that $T V_p\subset U_p$, and then by the above reasoning we may take $n\in N_{T} (y, V_p)\cap N_{T} (U, T^{- 1} V)$, and hence $n+ 1\in N_{T} (y, U_p)\cap N_{T} (U, V)$; it can be treated similarly
the other case that $U$ and $V$ can be viewed as opene subsets of $Y_c$ and $Y$, respectively.
 \end{proof}

Now we continue our proof. As $(Y, F)$ is topologically mixing, the system $(Y, F^2)$ is weakly mixing, and then by \cite[Proposition 3.8]{HKhKZh2015} we may choose $x^*\in Y$ such that $\omega_{\mathcal{N}_{F^2}} (x^*)= Y$, in particular, $\omega_{\mathcal{N}_{F^2}} (x^*)\supset Z$. Thus, by the above Claim, we obtain $Z\subset \omega_{\mathcal{N}_T} (x^*)$ and hence $Z= \omega_{\mathcal{N}_T} (x^*)$.
\end{proof}

Note that a dynamical system is
proximal if and only if it contains the unique fixed point, which is the only minimal point of the system \cite{AK}. Thus, as a direct corollary of Lemma \ref{Shar} and Theorem \ref{20151026}, we have:

\begin{cor} \label{0207}
There exists a non-proximal, non totally transitive, transitive compact system $(X, T)$ and a point $x_0\in X$ such that $\omega_{\mathcal{N}_T} (x_0)\neq
\omega_T (x)$ for all $x\in X$.
\end{cor}

Nevertheless is still open the following

\begin{ques B}
 Let $(X,T)$ be a weakly mixing system.   Is any  $\omega_{\mathcal{N}_T}(x)$ undecomposable into
 $\alpha$ disjoint closed, nonempty, positively $T$-invariant subsets, where $2\le \alpha\le \aleph_0$?
\end{ques B}

At the end of this section let us prove one more chaotical property of transitive compact systems in additional to already known
in
\cite{HKhKZh2015}.

Recall that a pair of points $x, y \in X$ is \emph{asymptotic} if $\lim_{n\rightarrow \infty} d (T^n x, T^n y)= 0$.
Denote by $\Asym_T (X)$ the set of all asymptotic pairs of points.
Any pair $(x, y)\in \Prox_T (X)\setminus \Asym_T (X)$ is called a \emph{Li-Yorke pair}. Recall that a dynamical system $(X, T)$ is \emph{Li-Yorke chaotic}
if there exists an uncountable set $S\subset X$ with $(S\times S)\setminus \Delta_2 (X)\subset \Prox_T (X)\setminus \Asym_T (X)$, where $\Delta_2 (X)= \{(x, x): x\in X\}$.

\begin{prop} \label{0200a}
Each transitive compact system $(X, T)$ is Li-Yorke chaotic.
\end{prop}

\begin{proof}
Clearly $(X, T)$ is transitive. Observe that we have assumed the state space to be not a singleton and in fact a compact metric space without isolated
points, then $(X, T)$ is a transitive system with $X$ infinite. Thus, the subset $\Asym_T (X)$ is a fist category subset of $X\times X$ by \cite[Corollary 2.2]{HuangYe2000}. It is easy to show that $\Prox_T (X)$ is a $G_\delta$ subset of $X\times X$, and applying \cite[Proposition 3.7]{HKhKZh2015} to the transitive compact system $(X, T)$ we have that $\Prox_T (x)$ is a dense subset of $X$ for each $x\in X$. Thus $\Prox_T (X)$ is a dense $G_\delta$ subset of $X\times X$, and then $\Prox_T (X)\setminus \Asym_T (X)$ is a second category subset of $X\times X$. Now applying the well-known Mycielski Theorem \cite[Theorem 1]{Myc64} we obtain an uncountable subset $S\subset X$ with $(S\times S)\setminus \Delta_2 (X)\subset \Prox_T (X)\setminus \Asym_T (X)$. That is, $(X, T)$ is Li-Yorke chaotic.
\end{proof}

\section{Weak transitive compactness and sensitivity for linear operators}

In this section we are considering the dynamics of linear operators on infinite dimensional spaces in relation to the properties studied
in previous sections. More precisely, we will show the equivalence of the topological weak mixing property with a weak version of transitive compactness. We obtain some results on transitive sensitivity too.

One should keep in mind that, for a linear dynamical system $(X, T)$, where $X$ is an infinite dimensional space, neither compactness nor
even local compactness of $X$ is satisfied. In particular, we are interested in the case where $X$ is an infinite dimensional separable
Banach space and $T:X\to X$ is a continuous linear map (in short, operator). In this framework, we will just write $(X,T)$ is an infinite
dimensional linear dynamical system. We recall that $X$ is a Banach space if it is a vector space endowed with a norm $\norm{\cdot}$ such that $X$ with the associated distance $d(x,y):=\norm{x-y}$ becomes a complete metric space. It is well known that $T: X\to X$ is an operator if and only if $\norm{T}:=\sup\{ \norm{Tx}: \norm{x}\leq 1\}<\infty$. We refer the reader to the books \cite{bm} and \cite{gepbook} for the theory of linear dynamics.

Note that all notations and concepts discussed in previous sections can be introduced into linear dynamics. We also introduce a weak version of dynamical compactness.
A linear system $(X,T)$ is called \emph{weakly dynamically compact with respect to the family
$\mathcal{F}$} if there exists a dense subset $X_0\subset X$ such that the $\omega_\mathcal{F}$-limit set
$\omega_{\mathcal{F}}(x)$ is nonempty for all $x\in X_0$. In particular,
 $(X, T)$ is called
\emph{weakly transitive compact}, if there exists a dense subset $X_0\subset X$ such that for any point $x\in X_0$ the $\omega_{\mathcal{N}_T}$-limit set
$\omega_{\mathcal{N}_T}(x)$ is nonempty,
in other words,
for any point $x\in X_0$ there exists
a point $z \in X$ such that $$ n_T(x,G) \cap N_T(U,V)\ne \varnothing$$
 for any  neighborhood
$G$ of $z$ and any opene subsets $U,V$ of $X$.

\begin{thm} \label{wmoptc}
Let $(X, T)$ be an infinite dimensional linear system. Then $(X, T)$ is weakly mixing if and only if it is weakly transitive compact.
\end{thm}

\begin{proof}
\textit{Sufficiency}. Suppose that $(X, T)$ is weakly transitive compact. Let $X_0\subset X$ be a dense subset such that, for each $x\in X_0$, there exists
 $z(x) \in X$ such that
 $$ n_T(x,G) \cap N_T(U,V)\ne \varnothing$$
 for any  neighborhood
$G$ of $z (x)$ and opene $U, V\subset X$. As $(X,T)$ is obviously transitive, by \cite[Theorem 5]{gepracsam} (see also \cite[Theorem 2.45]{gepbook}) to obtain the weak mixing property we just need to show that, for each opene $U\subset X$ and $0$-neighbourhood $W$, there is a continuous map $S: X\to X$ commuting with $T$ such that
\begin{equation} \label{(*)}
S(U)\cap W\ne \varnothing \ \ \ \mbox{ and } \ \ \ S(W)\cap U\neq \varnothing.
\end{equation}
Given an opene subset $U$ of $X$ and a $0$-neighborhood $W$, we fix $x\in U\cap X_0$ and $z(x)\in X$ accordingly to the weak transitive compactness of $(X,T)$. Since $0$-neighbourhoods are absorbing, we find a scalar $\lambda\ne 0$ such that $\lambda z (x)\in W$. Let $G$ be a neighbourhood of $z (x)$ such that $\lambda G\subset W$. By the hypothesis we can find
$$
m\in n_T(x,G) \cap N_T(\lambda W,U).
$$
That is, $T^mx\in G$ and so $\lambda T^mx\in W$; additionally, there exists $w\in W$ with $T^m\lambda w\in U$. Now pick $S:=\lambda T^m$, we have that $S$ commutes with $T$ and the property \eqref{(*)} is satisfied, therefore the system is weakly mixing.

\textit{Necessity}. Conversely, under the assumption of the weak mixing property for $(X,T)$, we know by \cite[Theorem 2.3]{bp99} (see also \cite[Theorem 3.15]{gepbook}) that there exists an increasing sequence $\{n_k: k\in \N\}\subset \N$ and a dense subset $X_0\subset X$ such that $T^{n_k}x\to 0$ for each $x\in X_0$ and, for arbitrary opene $U,V\subset X$, we can find $k\in \N$ such that $T^{n_k}(U)\cap V\ne\varnothing$. Thus, we obtain easily that $(X,T)$ is weakly transitive compact by selecting $z(x)=0$ for every $x\in X_0$.
\end{proof}

Concerning sensitivity, the situation is more complicated and, although we obtain some advances, three related problems are left open.

\begin{prop}
Let $(X, T)$ be an infinite dimensional linear, topologically transitive system. Then $(X, T)$ is thickly multi-sensitive, that is, there exists $\delta> 0$ such that $\bigcap_{i= 1}^k S_T (U_i, \delta)$ is thick for any finite collection of
opene $U_1, \dots, U_k\subset X$.
\end{prop}

\begin{proof}
Let $U_1, \dots, U_k$ be opene sets, and let $m\in \N$. Pick points $x_1, \dots, x_k$ such that $x_i \in U_i$ and choose $\varepsilon>0$ such
that $B_\varepsilon (x_i) \subset U_i$, where $B_\varepsilon (x_i)$ is the open ball of radius $\varepsilon$ centered at  $x_i$, for all $i \in \{1, \dots, k\}$. By a \emph{hypercyclic} vector we mean that its orbit is dense in the space $X$. Take a hypercyclic vector $u \in B_\varepsilon (0)$ by \cite[Theorem 2.19]{gepbook}, and
let $y_i=x_i+u$. Then $y_i \in U_i$ by the construction. Since $u$ is hypercyclic there is $n \in \mathbb{N}$, $n>m$, such that $\norm{T^n u}>(\norm{T}+1)^m$.
Then $\rho(T^{n-j}x_i, T^{n-j}y_i)=\norm{T^{n-j}(x_i-y_i)}=\norm{T^{n-j}u}>(\norm{T}+1)^{m-j}>1$ for all $i= 1, \dots, k$ and $j= 0, \dots, m- 1$. Hence $\{n, n- 1, \dots, n-m+1\} \subset  \bigcap_{i=1}^k S_T(U_i, 1)$,
and therefore $(X, T)$ is thickly multi-sensitive.
\end{proof}

\begin{prop} \label{tsop}
Let $(X,T)$ be an infinite dimensional linear system. Then the following conditions are equivalent:
\begin{enumerate}

\item For each $\delta> 0$, $(X, T)$ is transitively sensitive with a sensitive constant $\delta$.

\item There exists $\delta_0> 0$ such that $(X, T)$ is transitively sensitive with a sensitive constant $\delta_0$.

\item There exists $\delta_0> 0$ such that $S_T (W_0, \delta_0)\cap N_T (U, V)\neq \varnothing$ for any opene subsets $U, V$ of $X$ and any $0$-neighbourhood $W_0$.
\end{enumerate}
\end{prop}

\begin{proof}
We just need to show $(3) \Rightarrow (1)$. Indeed, let $\delta>0$ be arbitrary, and fix arbitrary opene $U, V, W$ of $X$. We select $\varepsilon>0$ and $x\in W$ such that $x+B_\varepsilon(0)\subset W$. Observing $S_T (\lambda W_0, \lambda \delta_0)$ for any scalar $\lambda\ne 0$, and so without loss of generality we assume $\delta>\delta_0$. Let $0<\varepsilon '< \frac{\delta_0 \varepsilon}{\delta}$, and set $W_0=B_{\varepsilon'}(0)$. By the hypothesis there are $y,z\in W_0$ and $n\in N_T(U,V)$ such that $\norm{T^ny-T^nz}> \delta_0$. Set $y'=x+\frac{\delta}{\delta_0} y$ and $z'=x+\frac{\delta}{\delta_0} z$. We have $y', z'\in W$ and $\norm{T^ny'-T^nz'}> \delta$. As opene $U, V, W\subset X$ are arbitrary, $(X, T)$ is transitively sensitive with a sensitive constant $\delta$.
\end{proof}

In this framework the weak mixing property implies transitive sensitivity too. The following result establishes a very close connection of transitivity with transitive sensitivity. We do not know, however, whether every transitive linear system is transitively sensitive.

\begin{prop} \label{trtsop}
Let $(X,T)$ be an infinite dimensional linear, topologically transitive system. If $(X, T)$ is not transitively sensitive, then there exists a dense open subset $U_0\subset X$ such that every $x\in U_0$ has a dense orbit.
\end{prop}

\begin{proof}
If $(X,T)$ is not transitively sensitive, by Proposition \ref{tsop} we find opene $U, V$ of $X$ and $\delta>1$ such that $\norm{T^nx}\leq \delta$ whenever $n\in N_T(U,V)$ and $\norm{x}\leq 1$.  We fix an arbitrary opene $V'\subset V$ and select an opene $\widehat{V}\subset V'$ and $\varepsilon>0$ such that $\widehat{V}+B_\varepsilon(0)\subset V'$.  Given $u\in U$, there is $\varepsilon'< \frac{\varepsilon}{\delta}$ such that $U':=u+B_{\varepsilon'}(0)\subset U$.  Since $T$ is transitive, there exists $m\in N_T(U',\widehat{V})\subset N_T(U,V)$. That is, we find $u'=u+w\in U'$ with $\norm{w}<\varepsilon'$ and $T^mu'\in \widehat{V}$. By the assumption $\norm{T^mw}\leq \delta \varepsilon'<\varepsilon$. Therefore, $T^mu=T^mu'-T^mw\in \widehat{V}+B_\varepsilon (0)\subset V'$. Since $u\in U$ and opene $V'\subset V$ are arbitrary, we obtain that the orbit of every element in $U$ is somewhere dense, thus everywhere dense by transitivity of the system. Finally, the open set $U_0:=\bigcup_{n\in\N}T^{-n}(U)$ is dense,
and every element in $U_0$ has a dense
orbit.
\end{proof}

There are (very difficult) examples of linear systems $(X,T)$ such that every non-zero element has a dense orbit \cite{read}, but it seems to unknown whether every linear system that admits an open set of elements whose orbit is dense is so that every non-zero element has a dense orbit. It is also worthy to mention that there are (also rare) examples of transitive but not weakly mixing linear systems \cite{delarosaread} (see also \cite{bm}), but as far as we know there are no examples of transitive non-weakly mixing linear systems such that every non-zero element has a dense orbit.

Concerning weak disjointness, observe that for each separable Banach space the family of all opene subsets admits a countable base, and then it is a routine to show that Theorem~\ref{th0201} holds true within linear systems too. Note that the intersection of finitely many thickly syndetically sets is still thickly syndetic, and that an interesting property is that every \emph{topologically ergodic} linear system $(X, T)$ (i.e., each element of $\mathcal{N}_T$ is a syndetic set) satisfies that each element of $\mathcal{N}_T$ is actually a thickly syndetic set (see the exercises in \cite[Chapter 2]{gepbook}). Thus any finite family $(X_1,T_1), \dots ,(X_k,T_k)$ of topologically ergodic linear systems is weakly disjoint and, moreover, the product system $(X_1\times \dots\times X_k, T_1\times \dots\times T_k)$ is topologically ergodic.

\medskip

\bibliographystyle{amsplain}




\end{document}